\documentclass[preprint,12pt]{elsarticle}

\usepackage{algorithm}
\usepackage{algorithmicx}
\usepackage{algcompatible}
\usepackage{algpseudocode} 
\usepackage{amssymb,bm,amsmath}
\usepackage{svg}
\usepackage{booktabs}
\usepackage{multirow}
\usepackage{array} 
\usepackage{mwe}
\usepackage{graphbox} 
\usepackage{geometry}
\newcolumntype{M}[1]{>{\centering\arraybackslash}m{#1}}
\usepackage{subcaption}
\usepackage{comment}
\usepackage{array}
\usepackage{ragged2e}
\usepackage{float}
\usepackage{hyperref}
\usepackage{tikz}
\usepackage{xcolor}
\usepackage{pgfplots}  
\pgfplotsset{compat=1.18}  
\definecolor{myblue}{RGB}{056,072,113}
\definecolor{myyellow}{RGB}{231,189,057}
\definecolor{myred}{RGB}{171,079,063}
\definecolor{mygreen}{RGB}{071,133,090}
\definecolor{mypurple}{RGB}{078,076,114}
\definecolor{mypink}{RGB}{211,137,161}

\newcommand{\mat}[1]{\mathbf{#1}}

\newcommand{\todo}[1]{\textcolor[RGB]{255, 0, 0}{[\textbf{TODO:} {\em #1}]}}

\journal{Journal of Computational Physics}

\begin{document}

\begin{frontmatter}



	\title{Matrix-Free Multigrid with Algebraically Consistent Coarsening on Adaptive Octrees}

	\author[inst1]{Mengdi Wang}
	\ead{mengdi.wang@gatech.edu}

	\author[inst1]{Yuchen Sun}
	\ead{yuchen.sun.eecs@gmail.com}

	\author[inst1]{Bo Zhu\corref{cor1}}
	\ead{bo.zhu@gatech.edu}

	\cortext[cor1]{Corresponding author.}

	\affiliation[inst1]{organization={School of Interactive Computing},
		addressline={Georgia Institute of Technology},
		city={Atlanta},
		postcode={30332},
		state={Georgia},
		country={USA}}

	\begin{abstract}
		We present a matrix-free GPU multigrid preconditioner with algebraically
		consistent coarsening for solving Poisson equations on adaptive octree grids
		with irregular domains. Within uniform-resolution regions, the coarsening
		satisfies the Galerkin principle. At T-junctions between refinement levels,
		we propose a flux-consistent coarse-grid correction that restores cross-level
		consistency while preserving the compact matrix-free representation. The
		coarse operators are stored in a compact matrix-free form suitable for
		parallel execution on GPUs. Numerical experiments demonstrate second-order
		accuracy, grid-independent convergence when used with PCG, and robust
		performance on cut-cell problems arising in fluid simulation. On a single
		NVIDIA RTX 4090 GPU, the solver achieves full-solve throughputs above 200
		million cells per second on analytical Poisson tests and above 70 million
		cells per second on pressure projection problems in fluid simulation.
	\end{abstract}




\end{frontmatter}


\section{Introduction}

In incompressible flow simulations based on the projection method
\cite{bell1989second, almgren1998conservative},
a pressure Poisson equation must be solved at every time step, and this solve
typically dominates the total computational cost \cite{bridson2015fluid}.
As simulations scale to domains with multiscale features, such as boundary layers,
sharp interfaces, or complex immersed solid geometries, the pressure solver
must remain efficient across many orders of magnitude in spatial resolution.
Elliptic equations of this type arise across a broad range of applications
beyond fluid dynamics, including electrostatics \cite{teunissen2018afivo,
	tomida2023athena++} and phase-change modeling \cite{zhang2019comparative}.

Adaptive Mesh Refinement (AMR) \cite{berger1984adaptive, berger1989local}
addresses this resolution challenge by concentrating computational resources
only where they are needed.
Among AMR approaches, octree-based adaptive grids have become a practical
choice for three-dimensional simulations
\cite{losasso2004simulating, popinet2003gerris, setaluri2014spgrid},
offering a natural hierarchical structure that supports efficient traversal,
refinement, and coarsening.
However, at refinement boundaries, cells of different sizes share faces,
creating \emph{T-junctions} that break the regularity assumed by standard
uniform-grid solvers \cite{theillard2013multigrid, weinzierl2018quasi}.

A further complication arises when solid objects occupy part of the domain.
In many applications, such as flow around airfoils or biological structures
\cite{saurabh2021scalable, hejazialhosseini2010high}, the Poisson equation
must be solved on an irregular domain.
A common approach is the \emph{cut-cell} method
\cite{johansen1998cartesian, Ng2009efficient, gibou2002second, zarifi2017positive},
which partitions cells crossed by the solid boundary into fluid and solid
portions. The result is a coefficient-weighted Poisson equation
\begin{equation}
	\nabla\cdot(\beta\nabla p) = b,
\end{equation}
where the coefficient $\beta$ is spatially varying and may be discontinuous
across the fluid-solid interface \cite{botto2013geometric, teunissen2019geometric}.
Alternatively, the immersed boundary method \cite{peskin2002immersed} introduces
spatially varying coefficients through a different mechanism.
In either case, the solver must handle coefficient discontinuities robustly
at every level of the multigrid hierarchy.

Multigrid methods are highly efficient for solving elliptic equations at
scale, owing to their optimal linear complexity \cite{brandt1977multi}.
While they can be used as standalone solvers, for symmetric positive-definite (SPD)
systems like the Poisson equation, they are often employed as highly effective
preconditioners for the preconditioned conjugate gradient (PCG) method, significantly
accelerating convergence.
Their integration with octree AMR has been studied extensively across scientific
computing and related fields.
\citet{popinet2003gerris, popinet2015quadtree} developed the Gerris and
Basilisk solvers, employing octree-based grids and geometric multigrid (GMG)
for Poisson equations in complex geometries.
\citet{sampath2010parallel} demonstrated parallel GMG on
octree meshes with second-order accuracy and good scalability for large-scale
simulations.
\citet{teunissen2018afivo} provides an open-source GMG solver for
adaptively refined quad/octrees, later extended \cite{teunissen2023geometric}
to handle non-graded meshes and irregular boundaries more accurately.
In astrophysical simulations, \citet{tomida2023athena++} incorporated
GMG into its AMR framework, and AMReX \cite{zhang2019amrex} provides scalable
block-structured AMR infrastructure with multigrid-based elliptic solvers for
a broad range of applications.
For graphics applications, octree-based Poisson discretizations with PCG
\cite{losasso2004simulating, aanjaneya2019efficient} and highly efficient
tile-based Poisson solvers such as SPGrid \cite{setaluri2014spgrid} have also been explored.

A recurring challenge across all these frameworks is the treatment of
T-junctions: standard V-cycles are not directly applicable, and full
approximation scheme (FAS) cycles together with flux-conservative
interpolation operators are needed to maintain conservation across refinement
boundaries \cite{theillard2013multigrid, weinzierl2018quasi, feng2018mass,
	teunissen2018afivo, tomida2023athena++}.
For immersed boundary and cut-cell problems, additional strategies include
boundary-local Fourier analysis to optimize relaxation parameters
\cite{coco2023ghost}, moving least-squares velocity transformations
\cite{ghomizad2021structured}, and geometry-aware prolongation operators
for immersed finite elements \cite{chu2024multigrid}.

A deeper issue, however, concerns the consistency of the coarse-level
operators themselves.
Geometric multigrid constructs coarse operators directly from the geometry.
In the interior of a continuous fluid domain, this construction naturally
satisfies the \emph{Galerkin principle} $\mat{A}^{l-1} = \mat{R}\mat{A}^l\mat{P}$
and is therefore equivalent to algebraic multigrid \cite{stuben2001review, shao2022fast}.
However, at boundaries where $\beta$ is discontinuous, such as in cut-cell
formulations, the geometric approach deviates from the Galerkin principle.
The resulting mismatch between fine and coarse operators introduces additional
residuals at these interfaces and can degrade convergence substantially.

Algebraic multigrid (AMG) \cite{stuben2001review, notay2010aggregation,
	yang2002boomeramg}
resolves this by constructing the coarse operator algebraically via the
Galerkin identity, ensuring consistency regardless of coefficient regularity.
This makes AMG a natural choice for cut-cell and variable-coefficient problems.
However, fully satisfying the Galerkin identity across all adaptive
refinement interfaces in a sparse matrix structure is computationally
expensive and memory-intensive, making it hostile to GPU execution
\cite{deriaz2023high, feng2018mass}. A practical compromise is therefore to
keep the coarsening algebraically consistent, in the sense of satisfying the
Galerkin principle, within uniform-resolution patches and at cut-cell
boundaries, and to replace the strict algebraic construction by a flux-consistent
coarse-grid correction at T-junctions.
Previous efforts to bring AMG to adaptive octree grids remain significantly
more challenging than GMG. \citet{weinzierl2018quasi} proposed a
quasi-matrix-free approach for octrees, but their formulation still requires
storing local $3\times 3$ or $5\times 5$ dense subsystems. On uniform grids,
\citet{shao2022fast} showed that storing one diagonal and three off-diagonal
coefficients per cell is sufficient to represent the Galerkin coarse operator
in a fully matrix-free form, but their approach does not extend to adaptive
refinement.

As problem sizes grow to hundreds of millions of cells, GPU acceleration is
essential for competitive solver performance.
GPU-accelerated multigrid and algebraic solvers have been explored both for
structured uniform grids \cite{antepara2024high, naumov2015amgx} and for
adaptive simulations \cite{lyu2024efficient, raateland2022dcgrid,
	bell2012exposing}, demonstrating that GPU-native implementations can
substantially reduce solve times relative to CPU-based methods.
These observations suggest that, for adaptive Poisson problems with irregular
domains, it is desirable to retain as much algebraic consistency in the
coarsening as possible, since such consistency is central to robust, nearly
grid-independent convergence, especially when multigrid is used as a
preconditioner for PCG. At the same time, a practical GPU implementation must
avoid the memory and access overhead of explicitly stored sparse matrices. This
motivates a matrix-free multigrid formulation, with algebraically consistent
coarsening, tailored to the adaptive octree setting.

Existing adaptive octree and AMR geometric multigrid solvers, such as
\citet{teunissen2018afivo}, the octree GMG library of \citet{teunissen2019geometric},
and \citet{tomida2023athena++}, construct coarse-grid operators directly from
the geometry. While effective on CPU architectures, these approaches do not
naturally extend to a compact, GPU-friendly matrix-free representation, and
their purely geometric coarsening deviates from the Galerkin principle at
cut-cell boundaries. These observations suggest a formulation that preserves
the efficiency of matrix-free structured operators while remaining robust on
adaptive cut-cell octrees. We therefore design a matrix-free multigrid
preconditioner with algebraically consistent coarsening for the Poisson
equation on adaptive octrees. Within uniform-resolution regions, we adopt the
four-coefficient matrix-free representation introduced by \citet{shao2022fast},
in which storing one diagonal and three off-diagonal coefficients per cell is
sufficient to represent the Galerkin coarse operator on uniform grids. At
adaptive T-junctions, where a strict algebraic construction is no longer
available in this compact form, we replace the standard coarse-grid correction
with a flux-consistent FAS-style correction that reuses the FAS restriction
pattern for our linear system. As shown in \autoref{sec:fas-multigrid}, a
standard V-cycle on adaptive octrees can double-count interface flux
contributions across refinement boundaries; the FAS-style correction resolves
this by incorporating the flux already generated during fine-level smoothing
into the coarse-grid right-hand side, restoring cross-level consistency.

The entire hierarchy is stored with only one diagonal and three face
coefficients per cell, avoiding any sparse-matrix assembly. A tile-based GPU
implementation maps naturally to the octree block structure, fusing stencil
evaluation, restriction, and residual computation into shared-memory kernels
that minimize global-memory traffic. The resulting preconditioner is used
inside PCG, which remains robust even on irregular cut-cell configurations
where a geometric multigrid preconditioner fails to converge.

Numerical experiments confirm that the solver achieves second-order spatial
accuracy and nearly grid-independent PCG convergence. On a single NVIDIA RTX 4090
GPU, the full-solve throughput reaches $211\sim 295$~M~cell/s on the analytical
sinusoidal Poisson tests of \autoref{subsec:sinusodial-poisson-test} and
$72\sim 97$~M~cell/s on the practical cut-cell pressure projection problems in
\autoref{sec:kinematic_flow}. Recent GPU projection solvers on uniform grids
already reach the low hundreds of millions of cells per second; for example,
\citet{lyu2024efficient} report roughly $235$~M~cell/s. Our results place
the proposed method in a similar throughput range while targeting substantially
more challenging adaptive-octree domains with irregular boundaries.

The remainder of the paper is organized as follows. \autoref{sec:methods}
introduces the matrix-free discretization and the multigrid formulation with
algebraically consistent coarsening. \autoref{sec:octree} describes the
treatment of adaptive octree T-junctions and the FAS-style coarse-grid
correction. \autoref{sec:results} reports accuracy, convergence, and
performance results on analytical Poisson, cut-cell projection, and
moving-obstacle flow problems. Finally, \autoref{sec:conclusion} discusses
limitations and future directions.

We have released the implementation based on the open-source fluid simulator
\textit{Cirrus} \cite{wang2025cirrus} at
\url{https://github.com/wang-mengdi/Cirrus-amg}.

\section{Matrix-Free Poisson Solver}
\label{sec:methods}

\subsection{Discretized Poisson's Equation}
\label{sec:matrix-free-laplacian}

When using the projection method \cite{bridson2015fluid} to solve the Navier-Stokes equations, we need to solve the pressure Poisson equation $-\nabla^2 p=b$, where $b=\nabla\cdot \bm u$ is the divergence of velocity $\bm u$. According to the divergence theorem, integrating the Poisson equation on a cell $\mathcal{C}$ is equivalent to calculating the flux integration over its boundary:
\begin{equation}
	\label{eq:laplacian-flux-integration}
	\int_{\mathcal{C}} -\nabla^2 p\,\mathrm{d}V = \oint_{\partial\mathcal{C}} -\nabla p \cdot \bm{n}\,\mathrm{d}S.
\end{equation}

In a MAC grid, $p$ and $b$ are stored at the cell centers. Assuming the cell size is $h$, the discrete form of the flux integration \eqref{eq:laplacian-flux-integration} becomes the sum of fluxes over its six faces, given by
\begin{equation}
	\label{eq:discretized-flux-integration}
	\int_{\mathcal{C}} -\nabla^2 p\,\mathrm{d}V \approx \sum_{j \in \mathcal{N}_i} \frac{p_i - p_j}{h} S_{i,j},
\end{equation}
where $\mathcal{N}_i$ represents the neighbors of cell $i$ and $S_{i,j}$ is the fluid area on the face between two cells. The gradient $\nabla p$ is calculated with central difference at face centers.

\autoref{eq:discretized-flux-integration} defines a symmetric seven-point linear system.
Such a system can be represented in a matrix-free form by storing, for each cell $i$,
a diagonal coefficient $c_i$ together with three off-diagonal coefficients
$c_{i,x-}$, $c_{i,y-}$, and $c_{i,z-}$ corresponding to its negative $x$-, $y$-, and $z$-faces. In other words,
\begin{equation}
	\label{eq:projection-poisson-system}
	\begin{cases}
		c_{i,x-}=-S_{i,x-}/h, \\
		c_{i,y-}=-S_{i,y-}/h, \\
		c_{i,z-}=-S_{i,z-}/h, \\
		c_i=-(c_{i,x-}+c_{i,y-}+c_{i,z-}+c_{i_1,x-}+c_{j_1,y-}+c_{k_1,z-}).
	\end{cases}
\end{equation}
Here $S_{i,x-}$ is the fluid area on the $x-$ face of cell $i$, and $i_1$, $j_1$, $k_1$ are its x+, y+, z+ neighbors.

In the ghost fluid method~\cite{enright2003robust}, each cell is classified into one of three types:
a fluid cell (inner cell), a Dirichlet cell (air, with boundary condition $p=0$), or a Neumann cell
(solid boundary, with boundary condition $\partial p/\partial \bm n=0$)~\cite{mcadams2010parallel}.
The coefficients $c$ introduced above are then assigned according to the cell type.

For a fluid cell, the off-diagonal coefficients are set to
$c_{i,x-}=c_{i,y-}=c_{i,z-}=-h$ and the diagonal entry to $c_i=6h$, since each face has area $h^2$.
At Neumann boundaries, the pressure difference $p_i - p_{i_1}$ vanishes in
\eqref{eq:discretized-flux-integration}, so the corresponding cross-term coefficient is set to zero,
e.g.\ $c_{i_1,x-}=0$.
At Dirichlet boundaries, we enforce $p_{i_1}=0$, which modifies the corresponding term in
\eqref{eq:discretized-flux-integration} to $-hp_i$.
Thus we set $c_{i_1,x-}=-h$, storing a nonzero cross-term even though the Dirichlet value is prescribed.

Since Neumann and Dirichlet cells are not degrees of freedom in the Poisson system, we set their diagonal entries $c_i=0$, indicating that the Laplacian operator
$-\nabla^2p$ is evaluated only on fluid cells.
The formulation of \eqref{eq:projection-poisson-system} naturally extends to the cut-cell case,
where solid and fluid may share portions of a cell face, and the coefficients are adjusted accordingly. It is also straightforward to verify that the seven-point linear system defined by
\eqref{eq:projection-poisson-system} is always positive semi-definite,
and becomes positive definite in the presence of Dirichlet boundaries.

\subsection{Preconditioned Conjugate Gradient}

We solve the linear system defined by \eqref{eq:projection-poisson-system} using the preconditioned conjugate gradient (PCG) method, as summarized in \autoref{alg:preconditioned-cg-solve}. In each iteration, the preconditioning operator \(\mathsf{Precond}(\bm{b})\) computes an approximate solution to the linear system \(\mat{A}\bm{x}=\bm{b}\). The simplest choice is the identity operator \(\mat{I}\), which reduces PCG to the standard conjugate gradient method. In this work, however, we use one multigrid \(\mu\)-cycle as the preconditioning operation, which will be described in the next subsection. 

The implementation of \autoref{alg:preconditioned-cg-solve} requires only three operations: the linear operator \(\mat{A} \bm{x}\), which in our problem corresponds to \(-\nabla^2\), the dot product, and pointwise operations. Therefore, PCG can be implemented using our matrix-free representation. Furthermore, there is a special case: if the matrix \(\mat{A}\) has a null space, such as when all boundaries in the computational domain are Neumann boundaries, we need to add a constant to the residual after calculating the residual in lines 4 and 11 of the algorithm, to ensure that the mean value of the residual vector is zero.

\begin{algorithm}
	\caption{Preconditioned Conjugate Gradient (PCG)}\label{alg:preconditioned-cg-solve}
	\begin{algorithmic}[1]
		\State \textbf{Input:} Matrix \(\mat{A}\), Right-hand side \(\bm{b}\)
		\State \textbf{Output:} Solution vector \(\bm{x}^{(k)}\)
		\State \(\bm{x}^{(0)} \gets 0\)
		\State \(\bm{r}^{(0)} \gets \bm{b} - \mat{A} \bm{x}^{(0)}\)
		\State \(\bm{z}^{(0)} \gets \mathsf{Precond}(\bm{r}^{(0)})\)
		\State \(\bm{p}^{(0)} \gets \bm{z}^{(0)}\)
		\State \(k \gets 0\)

		\While{$\|\bm{r}^{(k)}\| > \varepsilon_{\mathrm{tol}}$}
		\State \(\alpha_k \gets \frac{(\bm{r}^{(k)}, \bm{z}^{(k)})}{(\mat{A} \bm{p}^{(k)}, \bm{p}^{(k)})}\)
		\State \(\bm{x}^{(k+1)} \gets \bm{x}^{(k)} + \alpha_k \bm{p}^{(k)}\)
		\State \(\bm{r}^{(k+1)} \gets \bm{r}^{(k)} - \alpha_k \mat{A} \bm{p}^{(k)}\)
		\State \(\bm{z}^{(k+1)} \gets \mathsf{Precond}(\bm{r}^{(k+1)})\)

		\State \(\beta_k \gets \frac{(\bm{r}^{(k+1)}, \bm{z}^{(k+1)})}{(\bm{r}^{(k)}, \bm{z}^{(k)})}\)
		\State \(\bm{p}^{(k+1)} \gets \bm{z}^{(k+1)} + \beta_k \bm{p}^{(k)}\)

		\State \(k \gets k + 1\)
		\EndWhile
	\end{algorithmic}
\end{algorithm}

\subsection{Multigrid Preconditioner}
\label{sec:algebraic-multigrid}

We use a multigrid solver as the preconditioner for PCG.
Multigrid constructs a hierarchy of grids with progressively coarser resolutions by doubling the cell size at each level.
For example, if the finest grid is a \(256^3\) uniform grid, the multigrid hierarchy successively introduces coarser grids of size \(128^3\), \(64^3\), and so on.

The multigrid recursion used in this work is summarized in \autoref{alg:mu-cycle-multigrid}.
More generally, the multigrid procedure can be parameterized by a cycle index \(\mu\), which specifies how many recursive coarse-grid solves are performed at each level.
When \(\mu=1\), the scheme reduces to the classical \textit{V-cycle}, while \(\mu=2\) corresponds to the \textit{W-cycle}.
The basic idea of multigrid is to apply several smoothing operations on each level to eliminate high-frequency errors, then transfer the residual to the next coarser level and recursively solve the residual equation.
The coarse-grid correction is then prolongated back to the finer level to update the approximate solution.
Through this process, multigrid removes error components of different frequencies on different levels, leading to rapid convergence.

The superscript \( l \) in \autoref{alg:mu-cycle-multigrid} denotes the grid level, where a higher index corresponds to a finer grid. If the coarsest grid has a cell size of \( h_0 \), then the cell size at level \( l \) is \( h_0 / 2^l \).

The matrices $\mathbf{R}$ and $\mathbf{P}$ are \textit{restriction} and \textit{prolongation} matrices respectively. $\mathbf{R}$ transfers the linear system from the finer level $l$ to the coarser level $l-1$, while $\mathbf{P}$ transfers the solution from the coarser level $l-1$ to the finer level $l$. For example, if level $l$ has 1024 fluid cells and level $l-1$ has 128, then $\mathbf{R}$ is a $128 \times 1024$ wide matrix, and $\mathbf{P}$ is a $1024 \times 128$ tall matrix. We use a constant prolongation matrix \cite{shao2022fast}:
\begin{equation}
	P_{i,J} = \begin{cases}
		1, & \text{if } i \text{ is a child of } J, \\
		0, & \text{otherwise.}
	\end{cases}
\end{equation}
We denote fine-level cells with lowercase $i$ and coarse-level cells with uppercase $J$. This constant prolongation $\mat{P}$ is equivalent to directly copying parent cell values to their children.

On the other hand, the restriction matrix $\mat{R}$ is the transpose of $\mat{P}$ with a constant scaling coefficient $\alpha$:
\begin{equation}
	\mathbf{P} = \alpha \mathbf{R}^\top.
\end{equation}
Therefore all non-zero terms in $\mat{R}$ are $1/\alpha$.
The linear system $\mathbf{A}^l$ solved at each multigrid level is defined by the \textit{Galerkin Principle}:
\begin{equation}
	\label{eq:galerkin-principle}
	\mathbf{A}^{l-1} = \mathbf{R}^l \mathbf{A}^l \mat{P}^{l-1}.
\end{equation}

We employ the \textit{red-black Gauss-Seidel} (RBGS) method as our smoothing operator,
since it is well-suited for parallel execution and demonstrates good performance in our tests.
In our implementation, we apply $\nu_l=2$ RBGS iterations at the beginning and end of each level and apply $\nu_b=10$ RBGS iterations at the coarsest grid. To maintain operator symmetry, these two groups should be carried out in opposite orders (red-then-black vs. black-then-red). Similarly, the $\nu_b$ smoothing iterations at the coarsest level (line 2) are split into two phases with opposite color orders. One multigrid iteration performs a constant amount of work per level, therefore having linear complexity in the number of unknown variables.

The prolongation coefficient \(\beta\) allows for overshooting when transferring the solution of the residual equation to a finer level, further accelerating convergence. When using multigrid as a preconditioner for PCG, setting \(\beta = 2\) achieves the best convergence performance \cite{shao2022fast}. However, if multigrid is used as a standalone solver, \(\beta\) should be set to 1 to prevent oscillations.



\begin{algorithm}[H]
	\caption{$\mu$-Cycle Multigrid}
	\label{alg:mu-cycle-multigrid}
	\begin{algorithmic}[1]
		\State \textbf{Input:} Current level $l$, cycle index $\mu$, initial guess $\bm u^l$, right hand side (RHS) $\bm b^l$
		\State \textbf{Output:} Updated solution $\bm u^l$

		\If{$l=0$}
		\State $\bm u^0 \gets \mathsf{Smooth}(\nu_b, \bm u^0, \bm b^0)$
		\State \Return $\bm u^0$
		\EndIf

		\State $\bm u^l \gets \mathsf{Smooth}(\nu_l, \bm u^l, \bm b^l)$
		\State $\bm r^l \gets \bm b^l - \mat{A}^l \bm u^l$
		\State $\bm b^{l-1} \gets \mat{R}^l\bm r^l$

		\State $\bm u^{l-1}\gets \bm 0$
		\For{$k=1,\dots,\mu$}
		\State $\bm u^{l-1}\gets \mathsf{\mu\text{-}Cycle}(l-1, \mu, \bm u^{l-1}, \bm b^{l-1})$ \Comment{Recursive coarse solve}
		\EndFor

		\State $\bm u^{l}\gets \bm {u}^l+\beta \mat{P}^l\bm{u}^{l-1}$

		\State $\bm u^l \gets \mathsf{Smooth}(\nu_l, \bm u^l, \bm b^l)$

		\State \Return $\bm u^l$
	\end{algorithmic}
\end{algorithm}

\subsection{Algebraic and Geometric Multigrid}
\label{sec:algebraic-and-geometric-multigrid}

Under the Einstein summation convention, each element of $\mat{A}^{l-1}$ \eqref{eq:galerkin-principle} can be expressed as
\begin{equation}
	\label{eq:galerkin-mat-coarsen}
	\begin{aligned}
		A^{l-1}_{I,J} & =R_{I,i}A^l_{i,j}P_{j,J},                                            \\
		              & =\frac{1}{\alpha} P_{i,I}A^l_{i,j}P_{j,J}                            \\
		              & =\frac{1}{\alpha}\sum_{i\in \mathcal{C}_I,j\in\mathcal{C}_j}A_{i,j}.
	\end{aligned}
\end{equation}
Here $\mathcal{C}_I$ means the children of cell $I$ at finer level $l$.

If coarser cells $I,J$ are adjacent, there are four $(i,j)$ pairs with non-zero $A_{i,j}^l$ across the shared face in \eqref{eq:galerkin-mat-coarsen}. If all these cells are pure fluid in \eqref{eq:projection-poisson-system}, we have $A_{I,J}^{l-1}=4\cdot(-h)\cdot\frac{1}{\alpha}=-4h/\alpha$ with $h$ being the cell size at level $l$. Similarly, the diagonal entry $A_{I,I}^{l-1}$ is the sum of $8$ diagonal terms $6h$ and $24$ non-diagonal terms $-h$ (each face is counted twice), resulting in a total of $A_{I,I}^{l-1}=24h/\alpha$.

It's easy to see that $\alpha=2$ exactly satisfies the matrix entries directly defined by the integral form \eqref{eq:projection-poisson-system}, with the cell size being $2h$ at level $l-1$. In this case, the non-diagonal terms of $\mat{A}^{l-1}$ are $-2h$, and the diagonal terms are $12h$.
In Computer Graphics, a common practice is to ignore the cell size $h$ and set all off-diagonal entries to $-1$, corresponding to $\alpha=4$.

Therefore, in fluid-only systems, the matrix $\mat{A}^{l-1}$ given by the Galerkin principle \eqref{eq:galerkin-principle} is equivalent to the matrix given directly by the grid discretization \eqref{eq:laplacian-flux-integration}. The construction that obtains each level's linear system \( \mat{A}^{l-1} \) directly from the grid, without explicit reference to the fine-level operator, is referred to as \textit{geometric multigrid}. In contrast, the construction that enforces the Galerkin principle on the algebraic system at each level is referred to as \textit{algebraically consistent coarsening}, in the sense of algebraic multigrid \cite{stuben2001review}.

It can be observed that the geometric construction does not require explicitly storing the matrix entries \( c_i, c_{i,x-}, c_{i,y-}, c_{i,z-} \) in memory. Instead, they can be computed at runtime, making the geometric construction more efficient than an explicitly assembled algebraic multigrid implementation \cite{mcadams2010parallel}. However, if the computational domain contains non-fluid regions, the geometric construction no longer satisfies the Galerkin principle. For example, if child $i$ of cell \( I \) is non-fluid, the terms related to $i$ will vanish from \( A_{I,I}^{l-1} \), making it no longer \( 12h \). Enforcing algebraic consistency in the coarsening is therefore necessary to retain grid-independent convergence rates, even in the presence of irregular domains.

In this paper, we do not assemble any coarse-grid system explicitly. The coarsening is performed in a matrix-free form that nevertheless satisfies the Galerkin principle in uniform-resolution regions, and is combined with a flux-consistent coarse-grid correction at T-junctions (\autoref{sec:fas-multigrid}). Throughout the paper we refer to this construction as \emph{algebraically consistent coarsening}, in the sense of algebraic multigrid \cite{stuben2001review}, and to the overall preconditioner as a \emph{matrix-free multigrid preconditioner with algebraically consistent coarsening}.

\subsection{Matrix-Free Multigrid Coarsening}

\begin{algorithm}
	\caption{Coarsen}\label{alg:amg-coarsen}
	\begin{algorithmic}[1]
		\State \textbf{Input:} matrix-free entries of $\mat{A}^{l}$, coarse cell $(I,J,K)$
		\State \textbf{Output:} $c^{l-1}_{I,J,K}$,$c^{l-1}_{I,J,K,x-}$,$c^{l-1}_{I,J,K,y-}$,$c^{l-1}_{I,J,K,z-}$

		\State $c^{l-1}_{I,J,K} \gets 0$, $c^{l-1}_{I,J,K,x-} \gets 0$, $c^{l-1}_{I,J,K,y-} \gets 0$, $c^{l-1}_{I,J,K,z-} \gets 0$, $\text{active\_cnt} \gets 0$

		\For{$d_i,d_j,d_k \in \{0,1\}^3$}
		\State $(i,j,k) \gets (2I+d_i,2J+d_j,2K+d_k)$
		\If{$c^{l}_{i,j,k} \neq 0$}
		\State $\text{active\_cnt} \gets \text{active\_cnt} + 1$
		\State $c^{l-1}_{I,J,K}\gets c^{l-1}_{I,J,K}+\frac{1}{\alpha}c^l_{i,j,k}$ \Comment{diagonal}
		\EndIf

		\If{$d_i=1$}
		\If{$c^l_{i,j,k}\neq 0\land c^l_{i-1,j,k}\neq 0$}
		\State $c^{l-1}_{I,J,K}\gets c^{l-1}_{I,J,K}+\frac{2}{\alpha}c^l_{i,j,k,x-}$ \Comment{2 cross terms in diagonal entry}
		\EndIf
		\Else
		\State $c^{l-1}_{I,J,K,x-}\gets c^{l-1}_{I,J,K,x-}+\frac{1}{\alpha}c^{l}_{i,j,k,x-}$ \Comment{non-diagonal}
		\EndIf

		\If{$d_j=1$}
		\If{$c^l_{i,j,k}\neq 0\land c^l_{i,j-1,k}\neq 0$}
		\State $c^{l-1}_{I,J,K}\gets c^{l-1}_{I,J,K}+\frac{2}{\alpha}c^l_{i,j,k,y-}$ \Comment{2 cross terms in diagonal entry}
		\EndIf
		\Else
		\State $c^{l-1}_{I,J,K,y-}\gets c^{l-1}_{I,J,K,y-}+\frac{1}{\alpha}c^{l}_{i,j,k,y-}$\Comment{non-diagonal}
		\EndIf

		\If{$d_k=1$}
		\If{$c^l_{i,j,k}\neq 0\land c^l_{i,j,k-1}\neq 0$}
		\State $c^{l-1}_{I,J,K}\gets c^{l-1}_{I,J,K}+\frac{2}{\alpha}c^l_{i,j,k,z-}$ \Comment{2 cross terms in diagonal entry}
		\EndIf
		\Else
		\State $c^{l-1}_{I,J,K,z-}\gets c^{l-1}_{I,J,K,z-}+\frac{1}{\alpha}c^{l}_{i,j,k,z-}$\Comment{non-diagonal}
		\EndIf
		\EndFor

		\If{$\text{active\_cnt}=0$}
		\State $c^{l-1}_{I,J,K}\gets 0$
		\EndIf

	\end{algorithmic}
\end{algorithm}

According to \eqref{eq:galerkin-mat-coarsen}, $\mat{A}^{l-1}$ preserves symmetry and follows a seven-point stencil pattern, since $A_{I,J}^{l-1}$ can only be non-zero when $I=J$ or when $I$ and $J$ are adjacent. Therefore, the linear system at each level solved by the multigrid can be fully represented by our matrix-free method. The implementations of \( \mat{R} \) and \( \mat{P} \) are straightforward, so the multigrid method can be implemented in a matrix-free manner.



We refer to the process of calculating the matrix-free coefficients \( c^{l-1} \) from \( c^l \) as \textit{Coarsening}, as shown in \autoref{alg:amg-coarsen}. Note that we use the non-zero diagonal entry \( c_{i,j,k}^l\neq 0 \) to mark whether a cell is an unknown variable in the equation, or "active", which is true if and only if it contains fluid.

\section{Multigrid on Octree}
\label{sec:octree}


\begin{figure}[ht]
	\centering
	\includegraphics[width=0.45\textwidth]{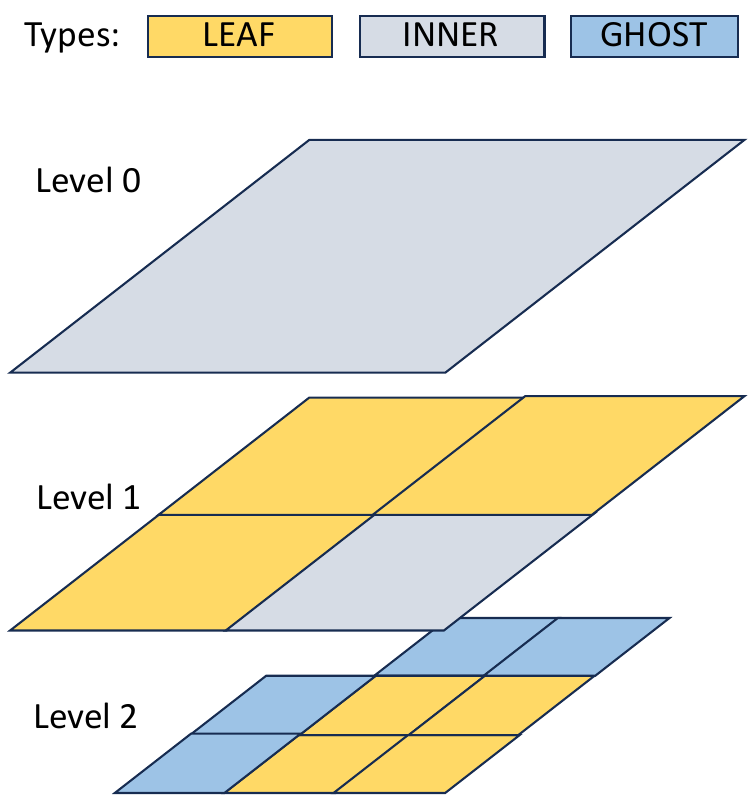}
	\caption{An 2D example of the octree structure we use, where different colors represent different types of cells. Actual octree grids are in 3D space.}
	\label{fig:octree-structure}
\end{figure}

Octree is a widely used data structure in computational physics for its ability to dynamically adjust grid resolution across different regions, thereby concentrating computational resources on critical areas with rich physics details. \autoref{fig:octree-structure} shows a 2D example of the adaptive octree we used in this paper.

As shown in the figure, we have three types of cells: leaf, inner, and ghost. A \textbf{leaf cell} represents a leaf node in the octree that stores active physical quantities, while an \textbf{inner cell} is the parent of 8 leaf cells or 8 inner cells. We enforce a level restriction that two face-neighboring leaf cells may differ by at most one level. When this case happens, a \textbf{ghost cell} is created next to the finer-level leaf cell to temporarily store the coarse leaf data, reducing expensive cross-level access when calculating the Laplacian operator. Next, we introduce how to define the Laplacian operator and the multigrid method on such an adaptive octree.

\subsection{Laplacian Operator on Octree}

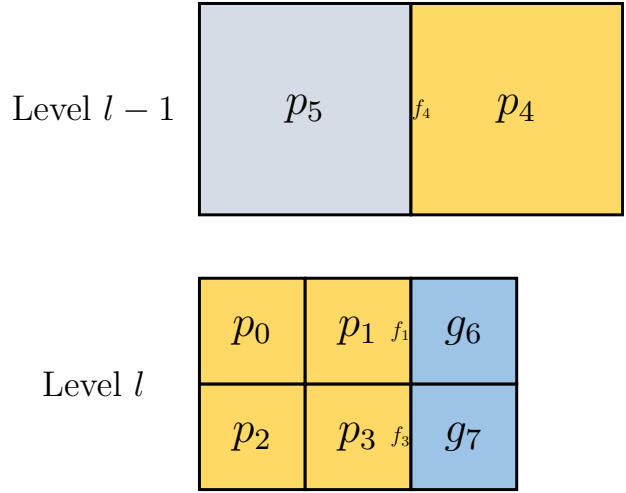
\begin{figure}[h]
	\centering
	\resizebox{0.6\textwidth}{!}{
		\begin{tikzpicture}
			\definecolor{innercolor}{rgb}{0.839, 0.863, 0.898}
			\definecolor{leafcolor}{rgb}{1.0, 0.851, 0.4}
			\definecolor{ghostcolor}{rgb}{0.616, 0.765, 0.902}
			\draw[black, fill=innercolor](0, 1.3) rectangle (1, 2.3);
			\draw[black, fill=leafcolor](1, 1.3) rectangle (2, 2.3);
			\draw[black, fill=leafcolor](0,0) rectangle (0.5, 0.5);
			\draw[black, fill=leafcolor](0.5,0) rectangle (1, 0.5);
			\draw[black, fill=leafcolor](0,0.5) rectangle (0.5, 1);
			\draw[black, fill=leafcolor](0.5,0.5) rectangle (1, 1);
			\draw[black, fill=ghostcolor](1,0) rectangle (1.5, 0.5);
			\draw[black, fill=ghostcolor](1,0.5) rectangle (1.5, 1);

			\node[scale=0.4] at (-0.5,0.5) {Level $l$};
			\node[scale=0.4] at (-0.5,1.8) {Level $l-1$};
			\node[scale=0.5] at (0.25,0.75) {$p_0$};
			\node[scale=0.5] at (0.75,0.75) {$p_1$};
			\node[scale=0.5] at (0.25,0.25) {$p_2$};
			\node[scale=0.5] at (0.75,0.25) {$p_3$};
			\node[scale=0.5] at (1.25,0.75) {$g_6$};
			\node[scale=0.5] at (1.25,0.25) {$g_7$};
			\node[scale=0.5] at (0.5,1.8) {$p_5$};
			\node[scale=0.5] at (1.5,1.8) {$p_4$};
			\node[scale=0.25] at (1.05,1.8) {$f_4$};
			\node[scale=0.25] at (0.95,0.75) {$f_1$};
			\node[scale=0.25] at (0.95,0.25) {$f_3$};

		\end{tikzpicture}
	}
	\caption{Illustration of a 2D T-junction. Cells $p_0\cdots p_3$ are at the finer level $l$ with cell size $h$. Cells $p_4$, $p_5$ are at the coarser level $l-1$ with cell size $2h$. $g_6$ and $g_7$ are ghost cells at level $l$ to avoid cross-level access. The fluxes of $\nabla p$ over $p_1-g_6$ boundary and $p_3-g_7$ boundary are $f_1$, $f_3$ respectively, and the flux of over $p_5-p_4$ boundary is $f_4$.}
	\label{fig:t-junction-2d}
\end{figure}

The main difference between an octree and a uniform grid is the presence of \textit{T-junctions}, as illustrated in the 2D example in \autoref{fig:t-junction-2d}. In the figure, leaf cells $p_0\cdots p_3$ are at the finer level $l$ with size $h$, leaf cell $p_4$ is at the coarser level $l-1$ with size $2h$, and inner cell $p_5$ is the parent of $p_0\cdots p_3$. Two ghost cells $g_6, g_7$ at level $l$ are next to $p_1$ and $p_3$ respectively.

The Laplacian operator requires defining the flux $f_1$ of \(\nabla p\) at the \( p_1 \)-\( g_6 \) boundary on level \( l \), similarly defining \( f_3 \) at the \( p_3 \)-\( g_7 \) boundary, and $f_4$ at the \( p_5 \)-\( p_4 \) boundary on level \( l-1 \).

With the matrix-free representation, we have
\begin{equation}
	\begin{aligned}
		f_1 & =c_{6,x-}(g_6-p_1), \\
		f_3 & =c_{7,x-}(g_7-p_3), \\
		f_4 & =c_{4,x-}(p_4-p_5),
	\end{aligned}
\end{equation}
A natural definition of $p_5$ is the average of its children, which means $p_5=(p_0+\cdots+p_3)/4$ in 2D, and finite difference \eqref{eq:projection-poisson-system} gives the coefficients
\begin{equation}
	\begin{aligned}
		c_{6,x-} & =\frac{S_{1}}{h},    \\
		c_{7,x-} & =\frac{S_{3}}{h},    \\
		c_{4,x-} & =\frac{S_1+S_3}{2h}.
	\end{aligned}
\end{equation}
Here $S_1$, $S_3$ are liquid areas over the two boundaries \( p_1 \)-\( g_6 \) and \( p_3 \)-\( g_7 \). In this paper, we continue to use the term \textit{area} to denote the length of 2D edges for brevity.

The Laplacian operator should conserve the flux on different levels by guaranteeing
\begin{equation}
	f_1+f_3=f_4,
\end{equation}
therefore
\begin{equation}
	\label{eq:flux-conservation}
	\frac{g_6-p_1}{h}S_1+\frac{g_7-p_3}{h}S_3=\frac{p_4-(p_0+\cdots+p_3)/4}{2h}(S_1+S_3)
\end{equation}
should hold for any values of $S_1$ and $S_3$ since we're solving Poisson's equation on an irregular domain. It forces us to take
\begin{equation}
	\label{eq:t-junction-ghost-values}
	\frac{g_6-p_1}{h}=\frac{g_7-p_3}{h}=\frac{p_4-(p_0+\cdots+p_3)/4}{2h},
\end{equation}
which means computing the finite difference $\nabla p$ at coarser level $l-1$ for the T-junction. It's easy to see that a similar version of \eqref{eq:t-junction-ghost-values} can also be applied to a 3D octree grid.

\subsection{Flux Deviation and FAS-Style Coarse-Grid Correction}
\label{sec:fas-multigrid}

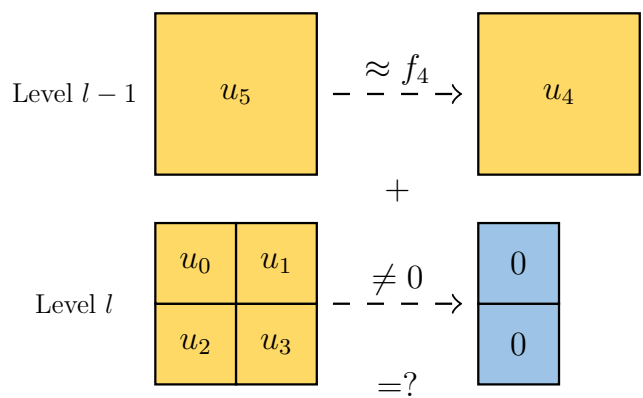
\begin{figure}[h]
	\centering
	\resizebox{0.6\textwidth}{!}{
		\begin{tikzpicture}
			\definecolor{innercolor}{rgb}{0.839, 0.863, 0.898}
			\definecolor{leafcolor}{rgb}{1.0, 0.851, 0.4}
			\definecolor{ghostcolor}{rgb}{0.616, 0.765, 0.902}
			\draw[black, fill=leafcolor](0,0) rectangle (0.5, 0.5); 
			\draw[black, fill=leafcolor](0.5,0) rectangle (1, 0.5); 
			\draw[black, fill=leafcolor](0,0.5) rectangle (0.5, 1); 
			\draw[black, fill=leafcolor](0.5,0.5) rectangle (1, 1); 
			\draw[black, fill=ghostcolor](2.0,0) rectangle (2.5, 0.5); 
			\draw[black, fill=ghostcolor](2.0,0.5) rectangle (2.5, 1); 
			\draw[dashed, ->] (1.1,0.5) -- (1.9,0.5);

			\node[scale=0.4] at (-0.5,0.5) {Level $l$};
			\node[scale=0.5] at (0.25,0.75) {$u_0$};
			\node[scale=0.5] at (0.75,0.75) {$u_1$};
			\node[scale=0.5] at (0.25,0.25) {$u_2$};
			\node[scale=0.5] at (0.75,0.25) {$u_3$};
			\node[scale=0.5] at (2.25,0.25) {$0$};
			\node[scale=0.5] at (2.25,0.75) {$0$};
			\node[scale=0.5] at (1.5,0.65) {$\neq 0$};
			\node[scale=0.5] at (1.5,1.2) {$+$};
			\node[scale=0.5] at (1.5,0.0) {$=?$};

			\begin{scope}[shift={(0,1.3)}]
				\draw[black, fill=leafcolor](0,0) rectangle (1, 1); 
				\draw[black, fill=leafcolor](2,0) rectangle (3, 1); 
				\draw[dashed, ->] (1.1,0.5) -- (1.9,0.5);

				\node[scale=0.4] at (-0.5,0.5) {Level $l-1$};
				\node[scale=0.5] at (0.5,0.5) {$u_5$};
				\node[scale=0.5] at (2.5,0.5) {$u_4$};
				\node[scale=0.5] at (1.5,0.65) {$\approx f_4$};
			\end{scope}
		\end{tikzpicture}
	}
	\caption{Flux deviation caused by directly applying the standard V-cycle at a T-junction on an adaptive grid.}
	\label{fig:mg-update-error}
\end{figure}

On adaptive grids with T-junctions, the coarse-grid operator is not constructed using the Galerkin formulation. As a result, consistency between grid levels is not automatically preserved, which may lead to flux deviation when a standard V-cycle is applied.

Consider applying a V-cycle at the T-junction shown in \autoref{fig:t-junction-2d}, which consists of three steps. First, perform RBGS smoothing on the finer level \( l \), as shown in the lower part of \autoref{fig:mg-update-error}. In this step, the value of cell $u_4$ remains zero, as propagated to its ghost children. It yields an approximate solution \( u_0, \dots, u_3 \), which produces a non-zero flux at the boundary between \( u_1, u_3 \), and \( u_4 \).

The second step of V-cycle restricts the residual on cells \( u_0, \dots, u_3 \) to level \( l-1 \) and smooths at level \( l-1 \), producing approximate solutions \( u_5 \) and \( u_4 \). Here, $u_5$ solves the residual equation with $r_5$ on the right-hand side, while $u_4$ solves the original equation $b_4$. Therefore, from the perspective of \( u_4 \), the boundary between \( u_5 \) and \( u_4 \) should approximately generate a correct flux \( f_4 \) as in the exact solution.

In the third step, \( u_5 \) is prolongated to level \( l \), which means adding \( u_5 \) to \( u_0, \dots, u_3 \) according to our definition of \( \mat{P} \). This effectively sums the flux contributions from both levels; therefore, the flux through the boundary between \( u_1, u_3 \), and \( u_4 \) deviates from the desired $f_4$, leading to slower convergence. In other words, the same interface flux is implicitly introduced once at the fine level through smoothing and once again at the coarse level through the multigrid correction, resulting in an overestimation of the flux across the T-junction interface.

The cause of this flux deviation is that cell \( u_4 \) is solved once as a boundary condition at level \( l \) and again as an unknown variable at level \( l-1 \), resulting in double the flux. To resolve this issue, we replace the original V-cycle by a coarse-grid correction that follows the restriction pattern of the Full Approximation Scheme (FAS), as described in \autoref{alg:fas-multigrid}. The FAS restriction pattern was originally proposed for nonlinear problems; we apply it here in a purely linear setting to restore cross-level flux consistency, and we therefore refer to this construction as an \emph{FAS-style} correction.

\begin{algorithm}[H]
	\caption{FAS-Style $\mu$-cycle Multigrid}
	\label{alg:fas-multigrid}
	\begin{algorithmic}[1]
		\State \textbf{Input:} Current level $l$, cycle parameter $\mu$, initial guess $\bm u^l$, right hand side (RHS) $\bm b^l$
		\State \textbf{Output:} Updated approximate solution $\bm u^l$

		\If{$l=0$}
		\State $\bm u^0 \gets \mathsf{Smooth}(\nu_b, \bm u^0, \bm b^0)$ \Comment{Or direct solve}
		\State \Return $\bm u^0$ \Comment{Coarsest grid solution}
		\EndIf

		\State $\bm u^l \gets \mathsf{Smooth}(\nu_l, \bm u^l, \bm b^l)$ \Comment{pre-smoothing}
		\State $\bm r^l \gets \bm b^l - \mat{A}^l \bm u^l$ \Comment{Residual}

		\State \textcolor{red}{$\bm u_\star^{l-1}\gets \mathsf{Avg}(\bm u^l)$} \Comment{Initial guess}

		\State $\bm b^{l-1} \gets \textcolor{red}{\beta}\mat{R}^l\bm r^l\textcolor{red}{+\mat{A}^{l-1} \bm u^{l-1}_\star}$
		\Comment{Restriction}

		\State $\bm u^{l-1} \gets \bm u_\star^{l-1}$

		\For{$k = 1,\dots,\mu$}
		\State $\bm u^{l-1}\gets \mathsf{FAS\text{-}Style\text{-}Multigrid}(l-1,\mu,\textcolor{red}{\bm u^{l-1}}, \bm b^{l-1})$ \Comment{Recursive solve}
		\EndFor

		\State $\bm u^{l}\gets \bm {u}^l+ \mat{P}^l\textcolor{red} {\left(\bm{u}^{l-1}-\bm{u}^{l-1}_\star\right)}$ \Comment{Prolongation}

		\State $\bm u^l \gets \mathsf{Smooth}(\nu_l, \bm u^l, \bm b^l)$ \Comment{post-smoothing}

		\State \Return $\bm u^l$

	\end{algorithmic}
\end{algorithm}

The differences with the original V-cycle are marked with red color in \autoref{alg:fas-multigrid}. Instead of solving the residual equation at level $l-1$, it solves the composite function $\beta\mat{R}^l\bm r^l+\mat{A}^{l-1} \bm u^{l-1}_\star$ that takes the T-junction flux term $\mat{A}^{l-1} \bm u^{l-1}_\star$ generated by the smoothing of level $l$ into account. Consequently, it prolongates the update $\bm{u}^{l-1}-\bm{u}^{l-1}_\star$ of the solution in order to remove the residual at level $l$. This modification ensures that the coarse-grid problem incorporates the flux contribution already generated on the fine grid, thereby preventing the double-counting of interface fluxes and restoring consistency across levels.

\subsection{Restriction and Prolongation Coefficients}

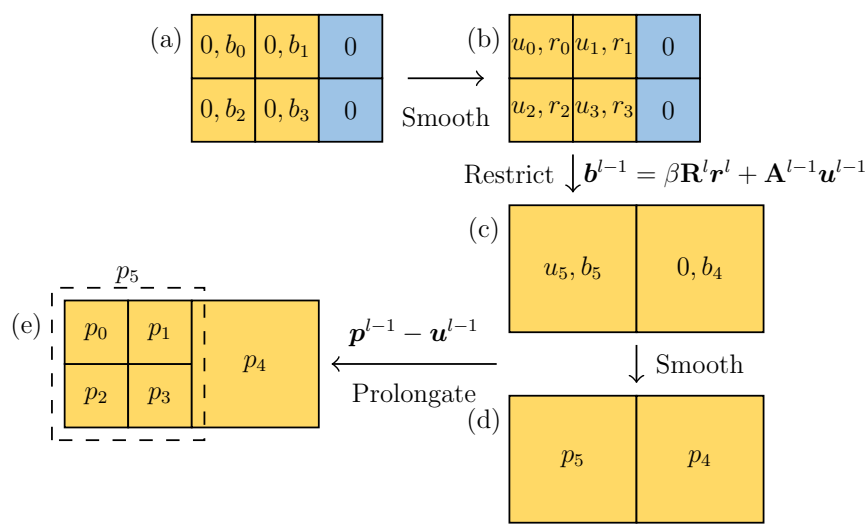
\begin{figure}[htpb]
	\centering
	\resizebox{0.8\textwidth}{!}{
		\begin{tikzpicture}
			\definecolor{innercolor}{rgb}{0.839, 0.863, 0.898}
			\definecolor{leafcolor}{rgb}{1.0, 0.851, 0.4}
			\definecolor{ghostcolor}{rgb}{0.616, 0.765, 0.902}

			\draw[black, fill=leafcolor](0,0) rectangle (0.5, 0.5); 
			\draw[black, fill=leafcolor](0.5,0) rectangle (1, 0.5); 
			\draw[black, fill=leafcolor](0,0.5) rectangle (0.5, 1); 
			\draw[black, fill=leafcolor](0.5,0.5) rectangle (1, 1); 
			\draw[black, fill=ghostcolor](1,0) rectangle (1.5, 0.5); 
			\draw[black, fill=ghostcolor](1,0.5) rectangle (1.5, 1); 

			\draw[->](1.7,0.5) -- (2.3,0.5);
			\node[scale=0.5] at (2,0.2) {Smooth};

			\node[scale=0.5] at (0.25,0.75) {$0,b_0$};
			\node[scale=0.5] at (0.75,0.75) {$0,b_1$};
			\node[scale=0.5] at (0.25,0.25) {$0,b_2$};
			\node[scale=0.5] at (0.75,0.25) {$0,b_3$};
			\node[scale=0.5] at (1.25,0.25) {$0$};
			\node[scale=0.5] at (1.25,0.75) {$0$};
			\node[scale=0.5] at (-0.2,0.8) {(a)};

			\begin{scope}[shift={(2.5,0)}]
				\draw[black, fill=leafcolor](0,0) rectangle (0.5, 0.5); 
				\draw[black, fill=leafcolor](0.5,0) rectangle (1, 0.5); 
				\draw[black, fill=leafcolor](0,0.5) rectangle (0.5, 1); 
				\draw[black, fill=leafcolor](0.5,0.5) rectangle (1, 1); 
				\draw[black, fill=ghostcolor](1,0) rectangle (1.5, 0.5); 
				\draw[black, fill=ghostcolor](1,0.5) rectangle (1.5, 1); 

				\draw[->](0.5,-0.1) -- (0.5,-0.4);
				\node[scale=0.5] at (0,-0.25) {Restrict};
				\node[scale=0.5] at (1.7,-0.25) {$\bm b^{l-1} = \beta\mat{R}^l\bm r^l+\mat{A}^{l-1} \bm u^{l-1}$};

				\node[scale=0.5] at (0.25,0.75) {$u_0,r_0$};
				\node[scale=0.5] at (0.75,0.75) {$u_1,r_1$};
				\node[scale=0.5] at (0.25,0.25) {$u_2,r_2$};
				\node[scale=0.5] at (0.75,0.25) {$u_3,r_3$};
				\node[scale=0.5] at (1.25,0.25) {$0$};
				\node[scale=0.5] at (1.25,0.75) {$0$};
				\node[scale=0.5] at (-0.2,0.8) {(b)};
			\end{scope}

			\begin{scope}[shift={(2.5,-1.5)}]
				\draw[black, fill=leafcolor](0,0) rectangle (1, 1); 
				\draw[black, fill=leafcolor](1,0) rectangle (2, 1); 

				\draw[->](1,-0.1) -- (1,-0.4);
				\node[scale=0.5] at (1.5,-0.25) {Smooth};

				\draw[->](-0.1,-0.25) -- (-1.4,-0.25);
				\node[scale=0.5] at (-0.75,-0.5) {Prolongate};
				\node[scale=0.5] at (-0.75,-0.0) {$\bm{p}^{l-1}-\bm{u}^{l-1}$};

				\node[scale=0.5] at (0.5,0.5) {$u_5,b_5$};
				\node[scale=0.5] at (1.5,0.5) {$0,b_4$};
				\node[scale=0.5] at (-0.2,0.8) {(c)};
			\end{scope}

			\begin{scope}[shift={(2.5,-3)}]
				\draw[black, fill=leafcolor](0,0) rectangle (1, 1); 
				\draw[black, fill=leafcolor](1,0) rectangle (2, 1); 

				\node[scale=0.5] at (0.5,0.5) {$p_5$};
				\node[scale=0.5] at (1.5,0.5) {$p_4$};
				\node[scale=0.5] at (-0.2,0.8) {(d)};
			\end{scope}

			\begin{scope}[shift={(-1.0,-2.25)}]
				\draw[black, fill=leafcolor](0,0) rectangle (0.5, 0.5); 
				\draw[black, fill=leafcolor](0.5,0) rectangle (1, 0.5); 
				\draw[black, fill=leafcolor](0,0.5) rectangle (0.5, 1); 
				\draw[black, fill=leafcolor](0.5,0.5) rectangle (1, 1); 
				\draw[black, fill=leafcolor](1,0) rectangle (2, 1); 
				\draw[black,dashed] (-0.1,-0.1) rectangle (1.1,1.1);

				\node[scale=0.5] at (0.25,0.75) {$p_0$};
				\node[scale=0.5] at (0.75,0.75) {$p_1$};
				\node[scale=0.5] at (0.25,0.25) {$p_2$};
				\node[scale=0.5] at (0.75,0.25) {$p_3$};
				\node[scale=0.5] at (1.5,0.5) {$p_4$};
				\node[scale=0.5] at (0.5,1.2) {$p_5$};
				\node[scale=0.5] at (-0.3,0.8) {(e)};
			\end{scope}

		\end{tikzpicture}
	}
	\caption{FAS-style multigrid update at a T-junction. (a) Initial state at finer level $l$. (b) Smooth on level $l$. Two ghost cells remain unchanged during the smoothing. (c) Restrict to coarser level $l-1$. The initial guess $u_5^\star$ is taken as the average of $u_0\dots u_3$, and the right-hand side is calculated as the sum of a residual term $\beta\mat{R}^l\bm r^l$ and a coarse flux term $\mat{A}^{l-1} \bm u^{l-1}_\star$. (d) Smooth on level $l-1$. (e) Prolongate the solution update $\bm{u}^{l-1}-\bm{u}^{l-1}_\star$ to level $l$ to get the final result $p_0\dots p_4$.}
	\label{fig:fas-update}
\end{figure}

The process of an FAS-style multigrid update is illustrated in \autoref{fig:fas-update}. We use $u$ to indicate the temporary solution and $p$ for the final result. The initial solution $u_5$ at level $l-1$ is taken as the average of four in our 2D example, which is $u_5=(u_0+\dots+u_3)/4$. As defined in \autoref{alg:fas-multigrid}, the prolongation update is $p_i=u_i+(p_5-u_5)$ for $i\in \{0,1,2,3\}$, therefore we have
\begin{equation}
	\label{eq:fas-flux-conservation}
	\frac{p_0+\dots +p_3}{4}=p_5.
\end{equation}
From the view of cell $p_4$, it means that it's solving the linear system with a correct final flux term with respect to $p_5-p_4$ at level $l-1$. Note that the prolongation (line 15) in \autoref{alg:fas-multigrid} no longer has the coefficient $\beta$, because any $\beta\neq 1$ in prolongation violates \eqref{eq:fas-flux-conservation}.

Fortunately, we can still achieve the same convergence as \(\beta = 2\) by applying \(\beta\) at the restriction step (\autoref{alg:fas-multigrid}, line 10). This is because the accuracy of the flux at the T-junction is guaranteed by the conservative conditions \eqref{eq:t-junction-ghost-values} and \eqref{eq:fas-flux-conservation}, and does not impose any requirements on the residual term solved for \( p_5 \). Since multigrid is a linear operator, multiplying the residual term by \( \beta = 2 \) is equivalent to multiplying the update step by \( \beta \), which has the same effect as using the prolongation coefficient.

Finally, according to \autoref{sec:algebraic-and-geometric-multigrid}, we use the restriction coefficient \(\alpha = 2\), because it corresponds exactly to the grid discretization in \eqref{eq:projection-poisson-system}. This means that when calculating \(\mat{A}^{l-1}\) using the Galerkin method, each non-zero element of the \( \mat{R} \) matrix is \( 0.5 \). This is especially important for octrees with T-junctions, as it ensures that, within the same level, the matrix-free coefficients obtained by coarsening the inner cells using \autoref{alg:amg-coarsen} remain consistent with those obtained directly from \eqref{eq:projection-poisson-system} for the leaf cells.

\subsection{GPU Implementation}

The solver is implemented on the GPU using the tile-based adaptive grid data structure introduced in our previous work~\cite{wang2025cirrus}.
The computational domain is divided into \(8\times8\times8\) \textit{tiles}, which serve as the smallest unit of memory allocation and kernel scheduling.
All cells inside a tile share the same type, which can be one of three categories: leaf, inner, or ghost.

For efficient stencil evaluation on the GPU, each CUDA thread block processes one tile.
Before computing the Laplacian operator, the required stencil data are loaded into shared memory, forming a \(10\times10\times10\) working region that includes the tile interior and one layer of neighboring cells.
The Laplacian operator is evaluated only on leaf and inner cells, while ghost cells appear only at the boundary of the shared-memory stencil.

Ghost cell values are not explicitly stored.
Instead, they are reconstructed on-the-fly according to \eqref{eq:t-junction-ghost-values}, which avoids cross-level memory accesses during stencil evaluation.
Because ghost values depend on the local refinement configuration, the same ghost cell may be reconstructed with different values when accessed by different neighboring tiles.

The discretized Poisson operator is implemented in a matrix-free manner.
Instead of assembling a global sparse matrix, each cell stores four coefficients corresponding to the diagonal term and the negative-direction face coefficients in the \(x\), \(y\), and \(z\) directions.
The positive-direction coefficients are obtained from neighboring cells during stencil evaluation.
This compact representation significantly reduces memory usage and allows the Laplacian operator to be evaluated on-the-fly during the solver iterations.

To further improve performance, multiple operations are fused into single GPU kernels whenever possible.
For example, each Red-Black Gauss-Seidel iteration computes the residual, diagonal term, and update step within a single kernel launch, minimizing redundant shared-memory loads.
Similarly, the residual computation and restriction step are fused into one kernel to reduce global memory traffic.

Finally, when the adaptive grid topology changes, we recompute the six same-level neighbors for each tile and cache their pointers directly inside the tile data structure.
During numerical computations these pointers are accessed directly, avoiding repeated hash-table queries and reducing irregular memory accesses.

\section{Results}
\label{sec:results}

\subsection{Overview of Test Cases and Grids}

\begin{figure}[!htbp]
	\centering
	\begin{minipage}{0.38\textwidth}
		\centering
		\includegraphics[width=\linewidth]{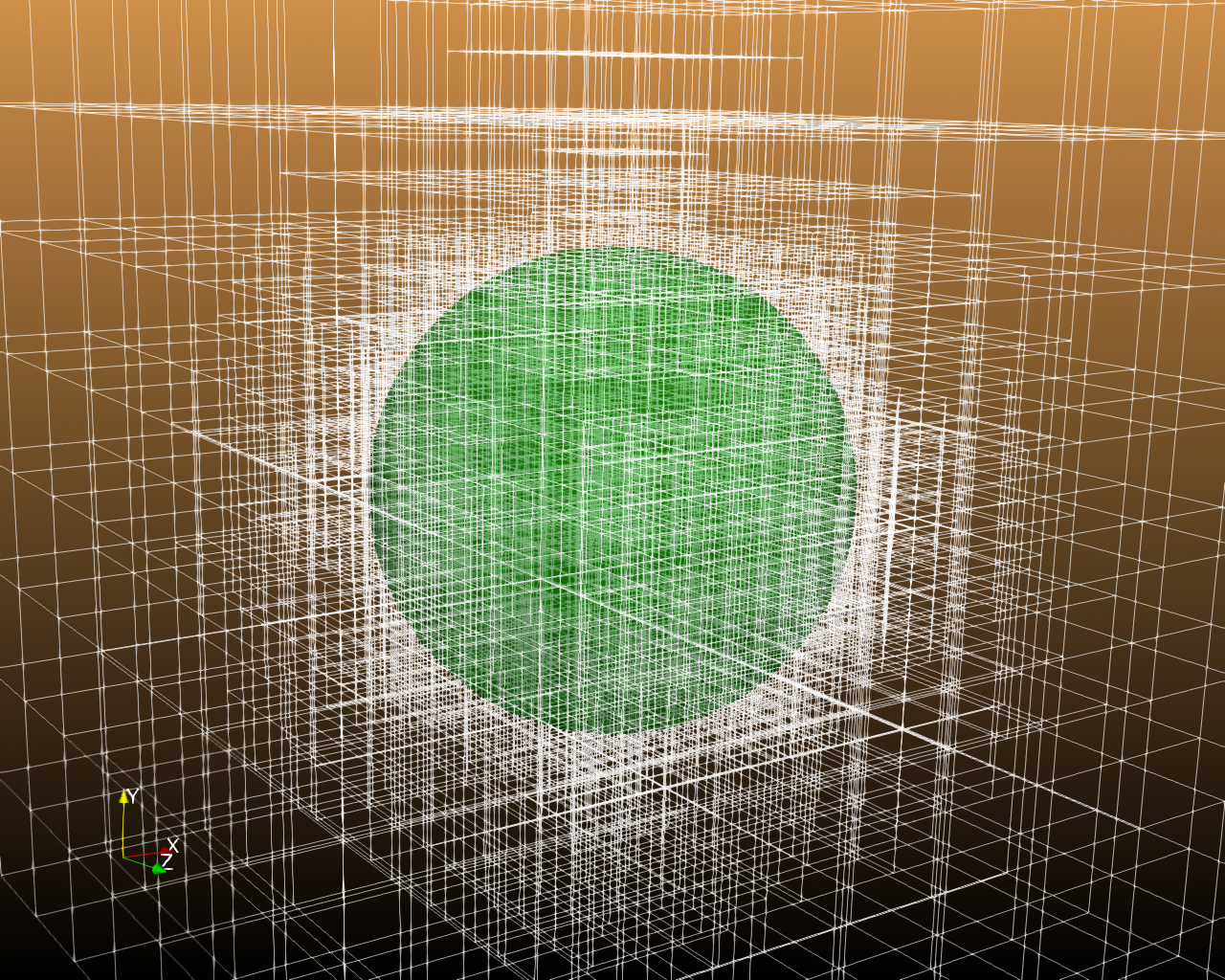}
	\end{minipage}\hspace{0.02\textwidth}
	\begin{minipage}{0.38\textwidth}
		\centering
		\includegraphics[width=\linewidth]{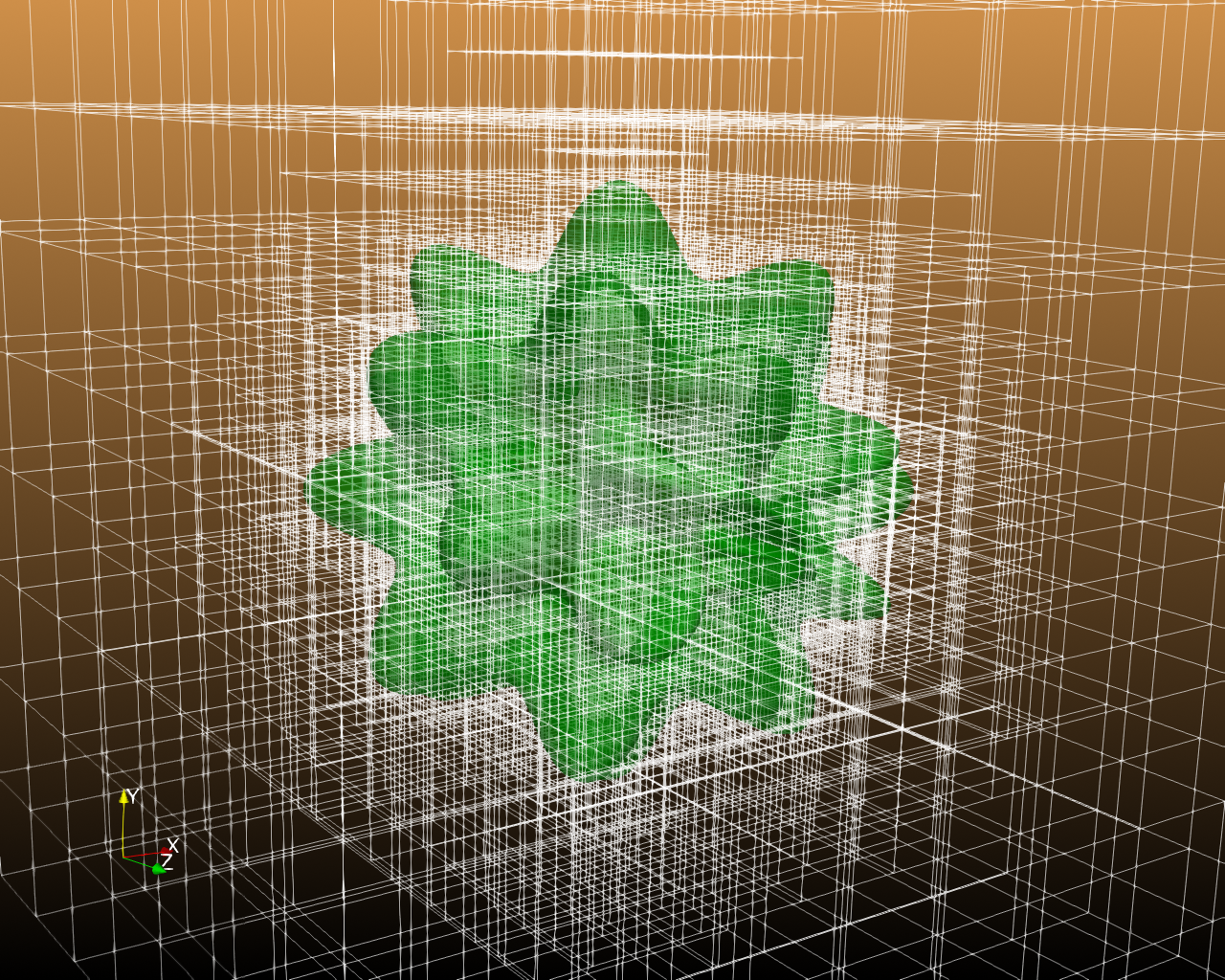}
	\end{minipage}
	\caption{
		Adaptive sphere grid (left) and star grid (right).
		Tile outlines are drawn in white and the isosurface $\phi=0$ is shown in green.
		Both configurations use $l_0=3$ with adaptive resolutions ranging from $64$ to $256$.
	}
	\label{fig:sphere_star_grids}
\end{figure}

\begin{figure}[!htbp]
	\centering
	\begin{minipage}{0.38\textwidth}
		\centering
		\includegraphics[width=\linewidth]{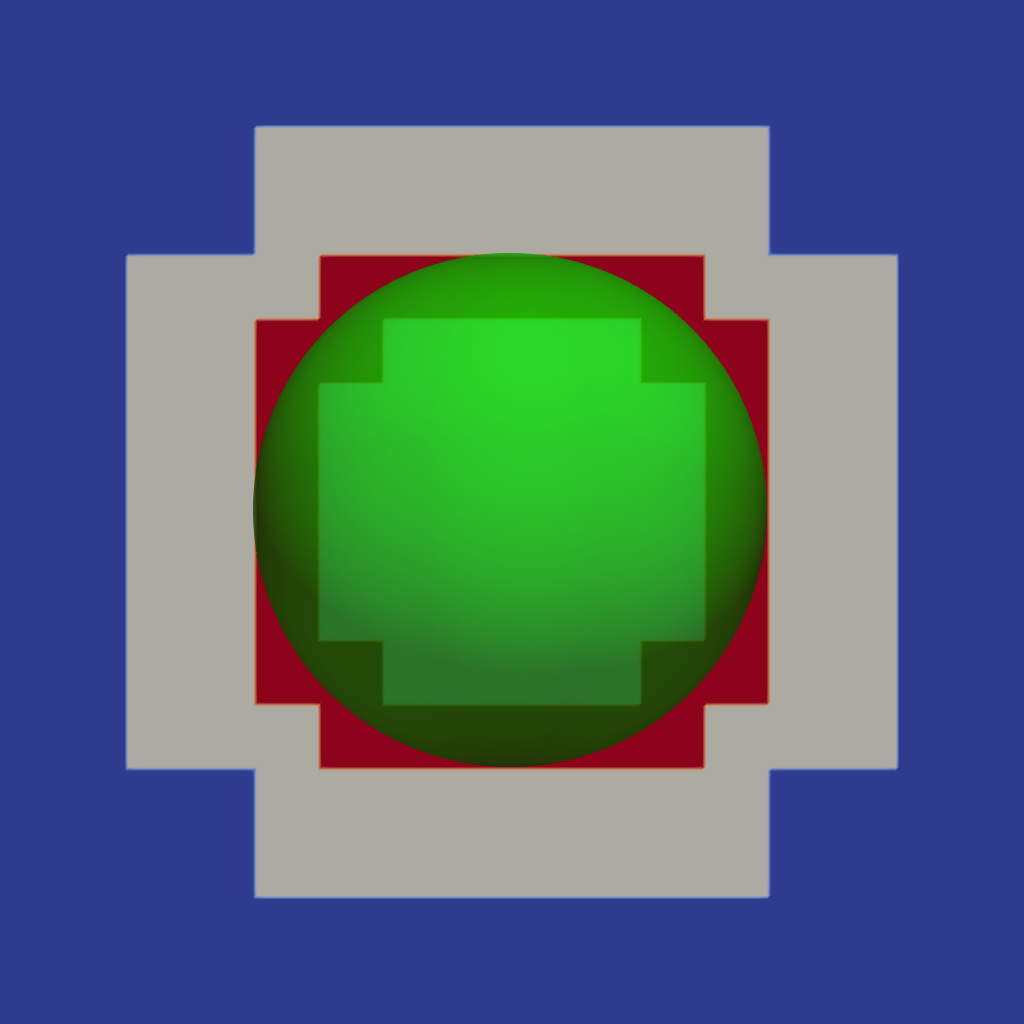}
	\end{minipage}\hspace{0.02\textwidth}
	\begin{minipage}{0.38\textwidth}
		\centering
		\includegraphics[width=\linewidth]{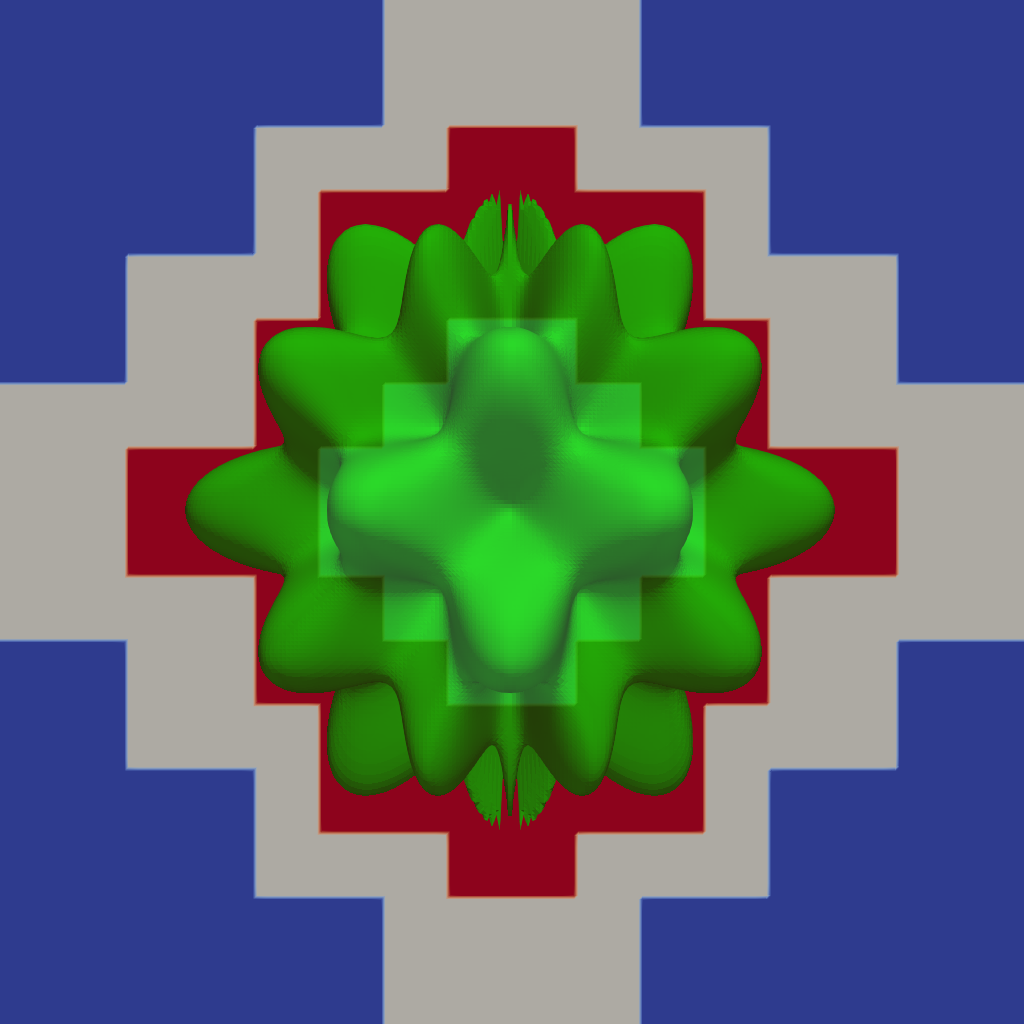}
	\end{minipage}
	\caption{
		Cross-sectional views at $x=0.5$ for the sphere grid (left) and the star grid (right) with $l_0=3$.
		Colors indicate refinement levels.
	}
	\label{fig:sphere_star_level_slices}
\end{figure}

In this section, we evaluate the proposed solver in terms of discretization accuracy, solver convergence, robustness on cut-cell domains, computational performance, and applicability in representative flow simulations.

All numerical tests are performed on the computational domain $[0,1]^3$.
We consider three types of grids: \textit{uniform}, \textit{sphere}, and \textit{star}.
The uniform grid serves as a baseline reference, while the sphere and star grids are designed to mimic a common situation in fluid simulation, where higher resolution is required only in a narrow band around solid boundaries.
In many applications, the physically important features, such as boundary layers and pressure variations induced by solid obstacles, are concentrated near the fluid-solid interface.
It is therefore natural to refine the grid only near the boundary and keep the rest of the domain coarse.
To reflect this usage pattern, we construct adaptive octree grids whose finest cells are concentrated in a narrow band around an implicit surface $\phi=0$.
In the sphere grid, this surface is a sphere, providing a simple and smooth boundary geometry, while in the star grid it is a star-shaped surface with stronger geometric variations.

The sphere shape is defined by the signed distance function
\begin{equation}
	\phi_\text{sphere}(\bm x)=|\bm{x}-(0.5,0.5,0.5)|-0.25,
\end{equation}
which describes a sphere centered at $(0.5,0.5,0.5)$ with radius $0.25$.
The star shape is defined by the signed distance function
\begin{equation}
	\phi_\text{star}(r, \theta, \phi) = r - \left[0.237 + 0.079 \cos(6\theta) \cos(6\phi)\right],
\end{equation}
also centered at $(0.5,0.5,0.5)$.

If a tile intersects the zero level set of the signed distance function, its target level is set to $l_0+2$. Otherwise, its target level is set to $l_0$. Since the maximum leaf level is two levels finer than the base level, the resulting grid contains three different leaf resolutions. We note that the solver itself is not limited to three refinement levels. The octree framework in our implementation naturally supports an arbitrary number of levels, and deeper refinement hierarchies can be obtained by specifying different target levels for the adaptive refinement.

Unless otherwise stated, we use base levels $l_0 \in \{2,3,4,5\}$.
\autoref{fig:sphere_star_grids} shows the sphere and star grids for $l_0=3$, where the narrow-band refinement around the implicit surface can be clearly observed.
\autoref{fig:sphere_star_level_slices} further illustrates the grid structure through cross-sectional views at $x=0.5$, where the colors indicate the refinement level of each tile.
In this configuration, the coarsest and finest leaf resolutions are $64$ and $256$, respectively.

All experiments were conducted on a workstation equipped with an NVIDIA RTX 4090 GPU (24GB memory) and an Intel i9-14900KF CPU.
The solver is implemented in CUDA and compiled using CUDA 12.9.
All numerical computations are performed in single precision, double precision is used only for the evaluation of dot products and volume-weighted RMS errors.
All experiments were executed on a single GPU.
The reported runtimes correspond to GPU wall-clock time.





\subsection{Accuracy of the Laplacian Operator}

\begin{figure}[!htbp]
	\centering
	\begin{tikzpicture}
		\begin{semilogyaxis}[
				width=0.6\textwidth,
				height=0.5\textwidth,
				xlabel={Root Cell Size},
				ylabel=RMS Error,
				xtick={1,2,3,4,5},
				xticklabels={
						1/16,
						1/32,
						1/64,
						1/128,
						1/256
					},
				cycle list={
						{color=myred, mark=*},
						{color=mygreen, mark=*},
						{color=myyellow, mark=*},
						{color=myblue, mark=*},
						{color=mypurple, mark=*},
						{color=mypink, mark=*}
					},
				yticklabel style={/pgf/number format/precision=0},
				legend pos=south west,
			]
			\addplot[color=myred, mark=*] coordinates {(2,7.50580227967458e-05) (3,1.214515125806391e-05) (4,1.7077968830644838e-06) (5,2.2220681128835373e-07)};
			\addlegendentry{sphere (2.80)};
			\addplot[domain=2:5, samples=100, dashed, thin, forget plot, color=myred] {2^(-8.012462337699981 + -2.803005046342034*x)};

			\addplot[color=myyellow, mark=*] coordinates {(2,3.854425509344119e-05) (3,1.1645500137404759e-05) (4,1.6592497916299014e-06) (5,2.1816860016535534e-07)};
			\addlegendentry{star (2.52)};
			\addplot[domain=2:5, samples=100, dashed, thin, forget plot, color=myyellow] {2^(-9.273436443409226 + -2.5205955279737906*x)};

			\addplot[color=myblue, mark=*] coordinates {(2,0.0001054703933407143) (3,1.3958991175765144e-05) (4,1.7944367734809868e-06) (5,2.27440762767143e-07)};
			\addlegendentry{uniform (2.95)};
			\addplot[domain=2:5, samples=100, dashed, thin, forget plot, color=myblue] {2^(-7.287995665593419 + -2.9530986239086827*x)};

		\end{semilogyaxis}
	\end{tikzpicture}
	\caption{The error between the Laplacian operator and the analytical Laplacian is measured with volume-weighted L2 norm on all cells. The convergence rates of grids are fitted with least squares (dashed lines) and shown in parentheses.}
	\label{fig:laplacian-error-vs-dx}
\end{figure}
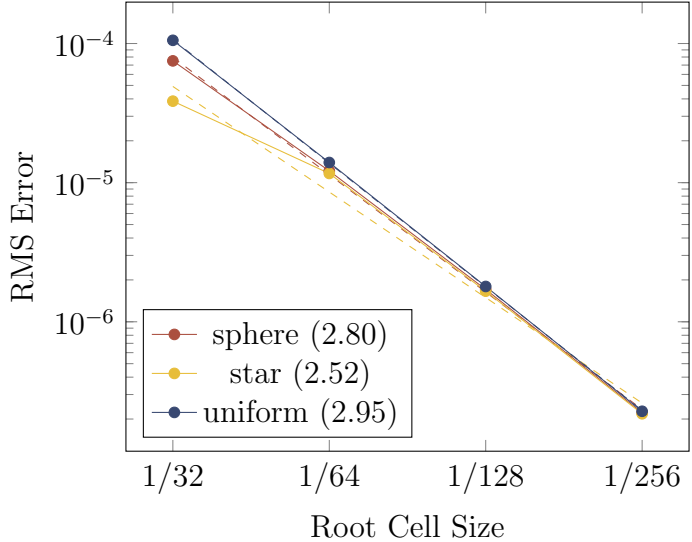

We first evaluate the discretization accuracy of the Laplacian operator on the adaptive octree grids.
To this end, we apply the numerical Laplacian to a smooth analytical function
\begin{equation}
	\begin{aligned}
		f(x,y,z)         & =(x^2 + y^2 + z^2) e^{-xyz},                                           \\
		\nabla^2f(x,y,z) & =2e^{-xyz}(x^2 + y^2 + z^2) - (xy + yz + zx)(x^2 + y^2 + z^2)e^{-xyz},
	\end{aligned}
\end{equation}
for which the exact Laplacian is known.

The numerical error is measured using the volume-weighted root mean square (RMS) value
\begin{equation}
	\mathrm{RMS}_V(e)=
	\sqrt{\frac{\sum_i V_i e_i^2}{\sum_i V_i}},
\end{equation}
where $e_i$ denotes the difference between the analytical and numerical Laplacian on cell $i$, and $V_i$ is the volume of that cell.
This metric accounts for the non-uniform cell sizes on adaptive grids and provides a consistent error measure across different grid configurations.

The Laplacian is evaluated on all leaf cells of the grid.
The outermost cells of the computational domain are treated as Neumann boundaries, while the remaining cells are considered interior cells. The discretization of the Laplacian operator and the treatment of boundary conditions follow the formulation introduced in \autoref{sec:matrix-free-laplacian}.

\autoref{fig:laplacian-error-vs-dx} shows the RMS error of the numerical Laplacian for different grid types and resolutions.
For all grid configurations, including uniform, sphere, and star grids, the method achieves better than second-order convergence with respect to the root cell size.

\subsection{Sinusoidal Poisson Solver Test}
\label{subsec:sinusodial-poisson-test}

\begin{figure}[htbp]
	\centering
	\includegraphics[width=0.38\linewidth]{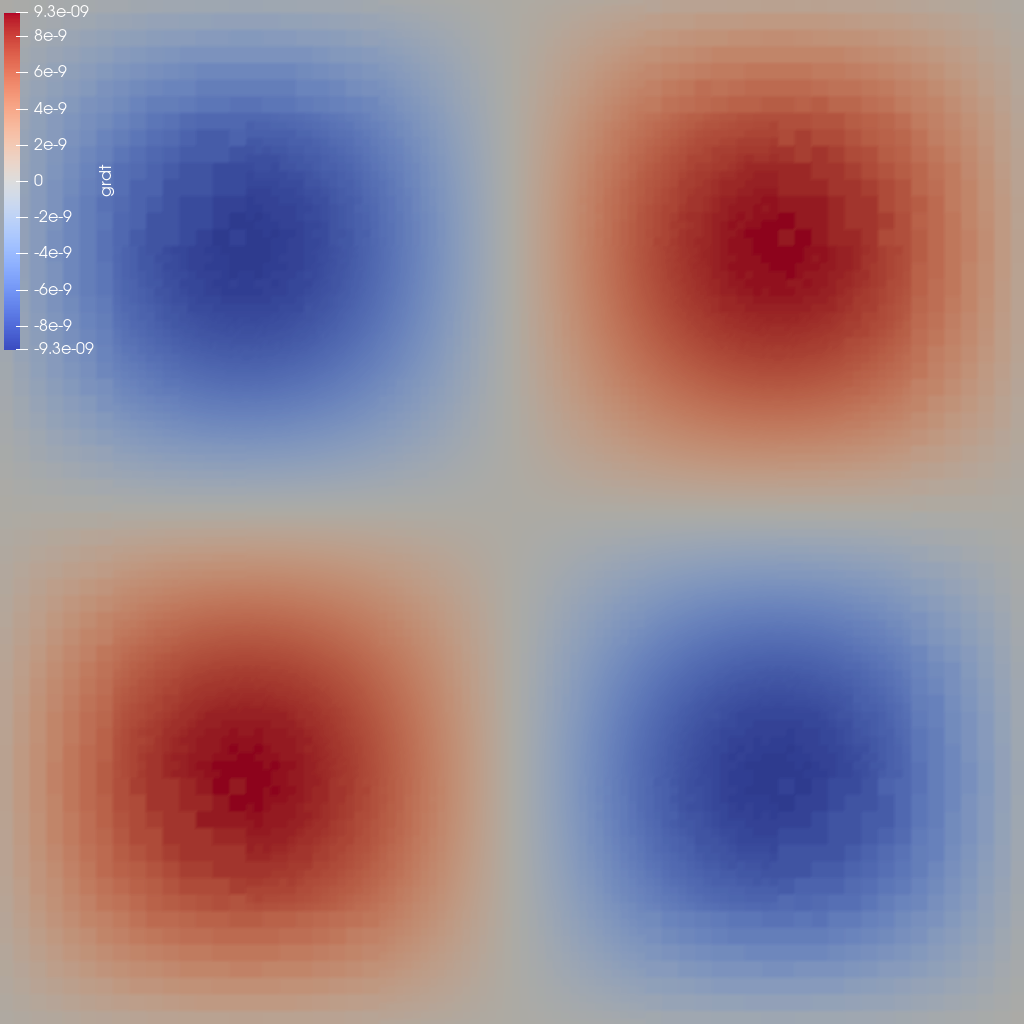}
	\hspace{0.02\textwidth}
	\includegraphics[width=0.38\linewidth]{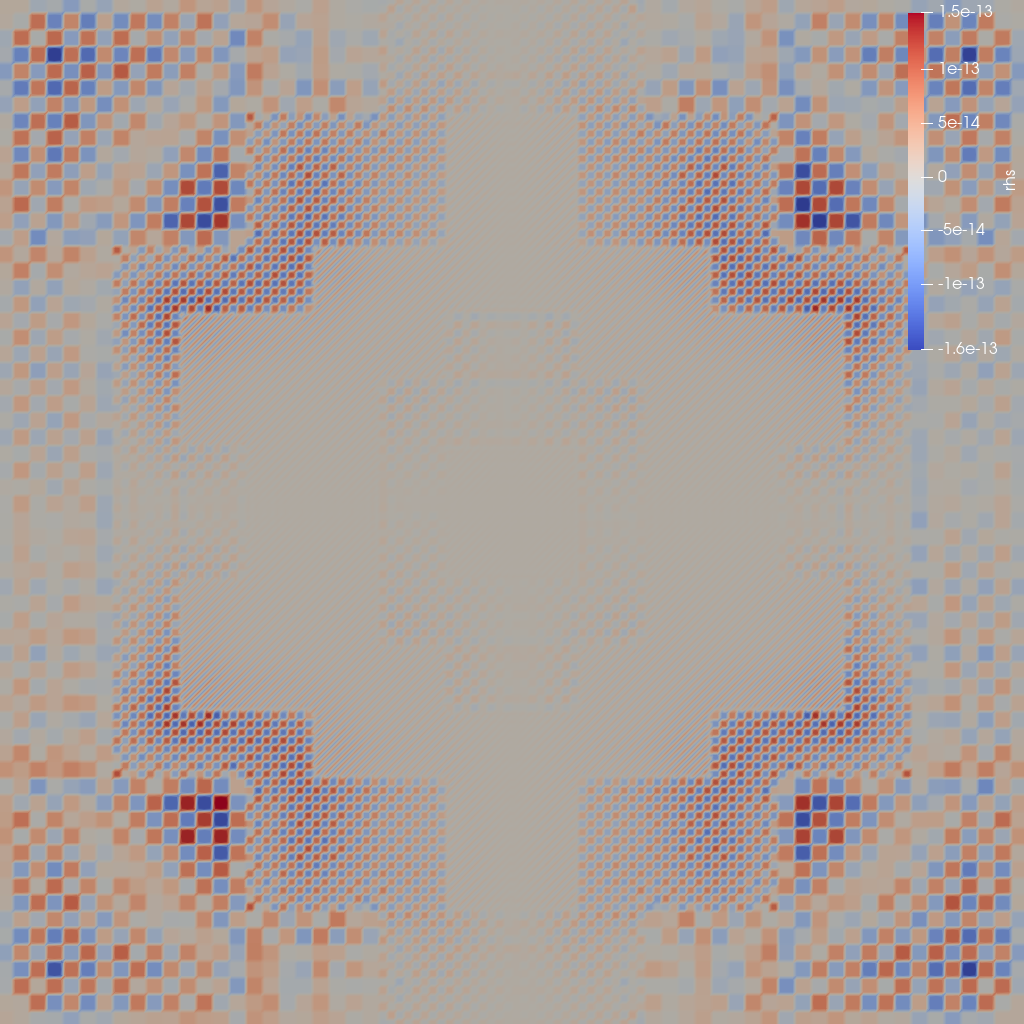}

	\caption{A slice at $z=0.5$ of the analytical solution $f(x,y,z)$ (left) and the numerical error (right) for the sinusoidal Poisson test. The grid is the star grid with $l_0=3$, having adaptive resolutions ranging from $64$ to $256$.}
	\label{fig:athena_star3_solution_residual}
\end{figure}

In this sinusoidal Poisson solver test, we evaluate the convergence behavior of the preconditioned conjugate gradient (PCG) solver using our matrix-free multigrid with algebraically consistent coarsening as a preconditioner. Solving Poisson equations is an important step in many fluid simulation algorithms, particularly in pressure projections. To assess the accuracy and convergence behavior of the solver, we consider the following Poisson equation with a sinusoidal analytical solution \cite{tomida2023athena++}:
\begin{equation}
	\nabla^2 f(x,y,z)=4 \pi G A \sin\left(\frac{2\pi x}{L_x}\right) \sin\left(\frac{2\pi y}{L_y}\right) \sin\left(\frac{2\pi z}{L_z}\right).
\end{equation}
The analytical solution to this equation is given by
\begin{equation}
	f(x,y,z)=- \frac{4 \pi G A}{\left( \frac{2\pi}{L_x} \right)^2 + \left( \frac{2\pi}{L_y} \right)^2 + \left( \frac{2\pi}{L_z} \right)^2} \sin\left(\frac{2\pi x}{L_x}\right) \sin\left(\frac{2\pi y}{L_y}\right) \sin\left(\frac{2\pi z}{L_z}\right) + \phi_0,
\end{equation}
where $L_x=L_y=L_z=1$, $G=A=1$, and $\phi_0=0$. We set the outermost layer of cells in the grid to Neumann boundary conditions and treat all other cells as interior cells.

The linear system is solved using PCG with our matrix-free multigrid preconditioner. We use red-black Gauss-Seidel relaxation as the smoother, with two pre-smoothing and two post-smoothing iterations per V-cycle level. The coarsest level uses ten smoothing iterations. A zero initial guess is used. For the accuracy and convergence tests in this subsection, the solver is run until the relative residual falls below $10^{-8}$, ensuring that the measured solution error reflects the discretization error rather than the solver residual.

A slice along $z=0.5$ of the analytical solution and the corresponding numerical error on the star grid are shown in \autoref{fig:athena_star3_solution_residual}. The error is computed as the difference between the numerical solution and the analytical solution. The checkerboard pattern visible in the error field is introduced by the red-black Gauss-Seidel relaxation used in the multigrid smoothing.

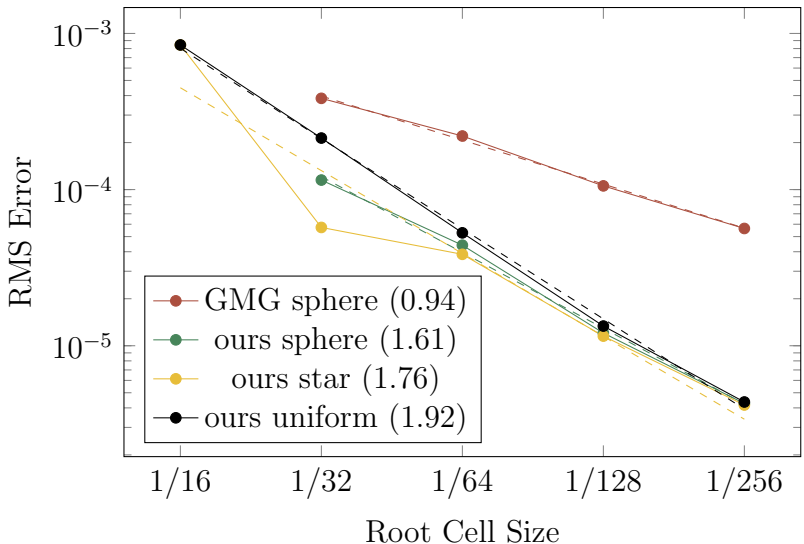
\begin{figure}[!htbp]
	\centering
	\begin{tikzpicture}
		\begin{semilogyaxis}[
				width=0.7\textwidth,
				height=0.5\textwidth,
				xlabel={Root Cell Size},
				ylabel=RMS Error,
				xtick={1,2,3,4,5},
				xticklabels={
						1/16,
						1/32,
						1/64,
						1/128,
						1/256
					},
				cycle list={
						{color=myred, mark=*},
						{color=mygreen, mark=*},
						{color=myyellow, mark=*},
						{color=myblue, mark=*},
						{color=mypurple, mark=*},
						{color=mypink, mark=*}
					},
				yticklabel style={/pgf/number format/precision=0},
				legend pos=south west,
			]

			\addplot[color=myred, mark=*] coordinates {(2,0.0003840788674912918) (3,0.00022035300933843978) (4,0.00010561748492014041) (5,5.6363469565433766e-05)};
			\addlegendentry{GMG sphere (0.94)};
			\addplot[domain=2:5, samples=100, dashed, thin, forget plot, color=myred] {2^(-9.426149787736634 + -0.9366678686633392*x)};


			\addplot[color=mygreen, mark=*] coordinates {(2,0.00011508008419653084) (3,4.4042708753048075e-05) (4,1.2063353417415805e-05) (5,4.267064749384272e-06)};
			\addlegendentry{ours sphere (1.61)};
			\addplot[domain=2:5, samples=100, dashed, thin, forget plot, color=mygreen] {2^(-9.788478220906502 + -1.6128023299641268*x)};

			\addplot[color=myyellow, mark=*] coordinates {(1,0.0008436836981357823) (2,5.722874652948698e-05) (3,3.8505026876251074e-05) (4,1.152763676862097e-05) (5,4.173389739398521e-06)};
			\addlegendentry{ours star (1.76)};
			\addplot[domain=1:5, samples=100, dashed, thin, forget plot, color=myyellow] {2^(-9.359583155727378 + -1.763032074534924*x)};

			\addplot[color=black, mark=*] coordinates {(1,0.0008436836123872862) (2,0.00021377620534548256) (3,5.2838098675036706e-05) (4,1.3369110737670753e-05) (5,4.363250377900705e-06)};
			\addlegendentry{ours uniform (1.92)};
			\addplot[domain=1:5, samples=100, dashed, thin, forget plot, color=black] {2^(-8.364686354486011 + -1.9189435881737404*x)};


		\end{semilogyaxis}
	\end{tikzpicture}
	\caption{Grid convergence for the sinusoidal Poisson solver test.
		The error between the numerical and analytical solutions is measured using the volume-weighted RMS error over all cells.
		GMG denotes the non-conservative geometric multigrid solver used for comparison.
		Dashed lines show least-squares fits used to estimate convergence rates (shown in parentheses).}
	\label{fig:athena-sin-grid-convergence}
\end{figure}

\autoref{fig:athena-sin-grid-convergence} shows the grid convergence behavior of the solver.
The error between the numerical and analytical solutions is measured using the volume-weighted RMS error over all cells.
For comparison, we also include results obtained with the non-conservative geometric multigrid method (GMG) \cite{losasso2004simulating}.
While GMG exhibits approximately first-order convergence on adaptive grids, our multigrid with algebraically consistent coarsening maintains convergence behavior close to second order across all grid configurations, demonstrating the benefit of conservative discretization across adaptive interfaces.

\begin{figure}[!htbp]
	\centering

	\begin{minipage}{0.49\linewidth}
		\centering
		\begin{tikzpicture}
			\begin{semilogyaxis}[
					width=\linewidth,
					height=0.72\linewidth,
					xlabel={PCG Iterations},
					ylabel={Residual},
					yticklabel style={/pgf/number format/precision=0},
					xtick={0,1,2,3,4,5,6},
					cycle list={
							{color=mygreen, mark=*},
							{color=myyellow, mark=*},
							{color=mypurple, mark=*}
						},
					legend to name=poissonpcglegend,
					legend columns=3,
				]


				\addplot[color=mygreen, mark=*] coordinates {(1,0.5496891) (2,0.04584512) (3,0.002557981) (4,0.0001483811) (5,5.537145e-06) (6,2.283441e-07)};
				\addplot[domain=1:6, samples=100, dashed, thin, forget plot, color=mygreen] {2^(3.887468664150532 + -4.261386162716564*x)};
				\addlegendentry{sphere $256\!\sim\!1024$}

				\addplot[color=myyellow, mark=*] coordinates {(1,0.2565234) (2,0.02135305) (3,0.001071922) (4,6.898856e-05) (5,2.815253e-06) (6,1.461549e-07)};
				\addplot[domain=1:6, samples=100, dashed, thin, forget plot, color=myyellow] {2^(2.5764469806163692 + -4.181147796085006*x)};
				\addlegendentry{star $64\!\sim\!256$}


				\addplot[color=mypurple, mark=*] coordinates {(1,0.5518096) (2,0.04576801) (3,0.00258618) (4,0.0001528315) (5,5.702078e-06) (6,2.376984e-07)};
				\addplot[domain=1:6, samples=100, dashed, thin, forget plot, color=mypurple] {2^(3.8721450238989856 + -4.249301029542801*x)};
				\addlegendentry{star $256\!\sim\!1024$}

			\end{semilogyaxis}
		\end{tikzpicture}
	\end{minipage}
	\hfill
	\begin{minipage}{0.49\linewidth}
		\centering
		\begin{tikzpicture}
			\begin{semilogyaxis}[
					width=\linewidth,
					height=0.72\linewidth,
					xlabel={PCG Iterations},
					ylabel={RMS Error},
					ytick pos=right,
					yticklabel pos=right,
					ylabel style={at={(axis description cs:1.22,0.5)},anchor=center},
					yticklabel style={/pgf/number format/precision=0},
					xtick={0,1,2,3,4,5,6},
					cycle list={
							{color=mygreen, mark=*},
							{color=myyellow, mark=*},
							{color=mypurple, mark=*}
						},
				]


				\addplot coordinates {
						(1,0.0030529857052860195) (2,0.00020304060766122942) (3,4.784492599994451e-06) (4,4.2151150898659714e-06) (5,4.264249348376803e-06) (6,4.264577629870643e-06)};

				\addplot coordinates {
						(1,0.0030324003567158584) (2,0.0001588338036648425) (3,3.630860097581442e-05) (4,3.852898977564876e-05) (5,3.851305679452261e-05) (6,3.8503981500021774e-05)};


				\addplot coordinates {
						(1,0.0030498592358039135) (2,0.00020125635089323737) (3,4.799541956799909e-06) (4,4.132355445195571e-06) (5,4.178178936129922e-06) (6,4.1780023152005805e-06)};

			\end{semilogyaxis}

		\end{tikzpicture}
	\end{minipage}

	\vspace{0.3em}

	\pgfplotslegendfromname{poissonpcglegend}

	\caption{PCG convergence for the sinusoidal Poisson solver test on representative adaptive grids.
		Left: PCG residual versus iterations.
		Right: RMS error between the numerical and analytical solutions versus PCG iterations.
		The solver maintains a grid-independent convergence rate with a residual reduction factor of approximately $1/20$ per iteration and reaches the discretization error level within about three iterations.
		The curves for sphere $256\!\sim\!1024$ and star $256\!\sim\!1024$ nearly overlap in the right plot and are therefore visually indistinguishable.}
	\label{fig:athena-sin-pcg-convergence}
\end{figure}
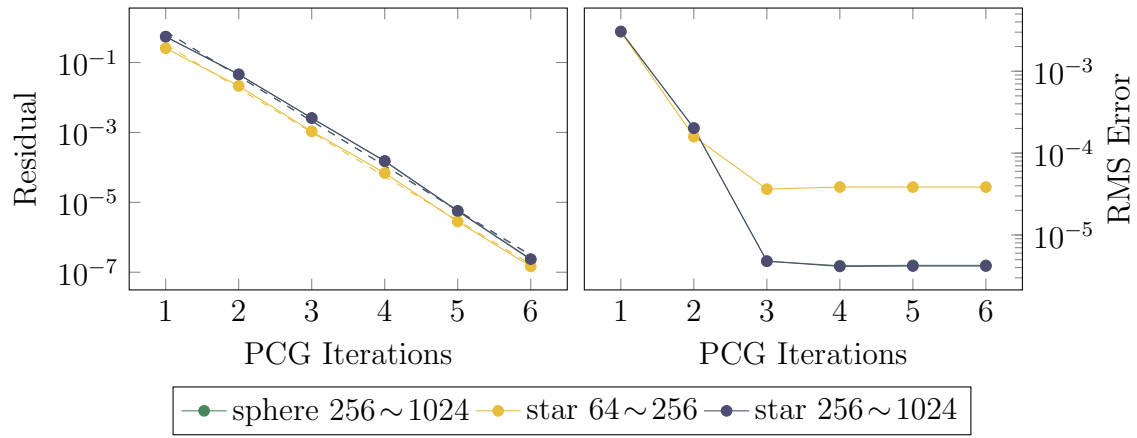

Finally, \autoref{fig:athena-sin-pcg-convergence} shows the convergence behavior of the PCG solver, terminated when the relative residual falls below $10^{-8}$.
Across all tested grids, the solver reaches the discretization error level within about three PCG iterations, while the residual decreases with a grid-independent reduction factor of approximately $1/20$ per iteration.
The fitted per-iteration reduction factors are $19.18$ for sphere $256\!\sim\!1024$, $18.14$ for star $64\!\sim\!256$, and $19.02$ for star $256\!\sim\!1024$.

\subsection{Static Pressure Projection Test with Cut Cells}
\label{subsec:static-projection-test}

\begin{figure}[htbp]
	\centering

	\begin{minipage}{0.38\textwidth}
		\centering
		\includegraphics[width=\linewidth]{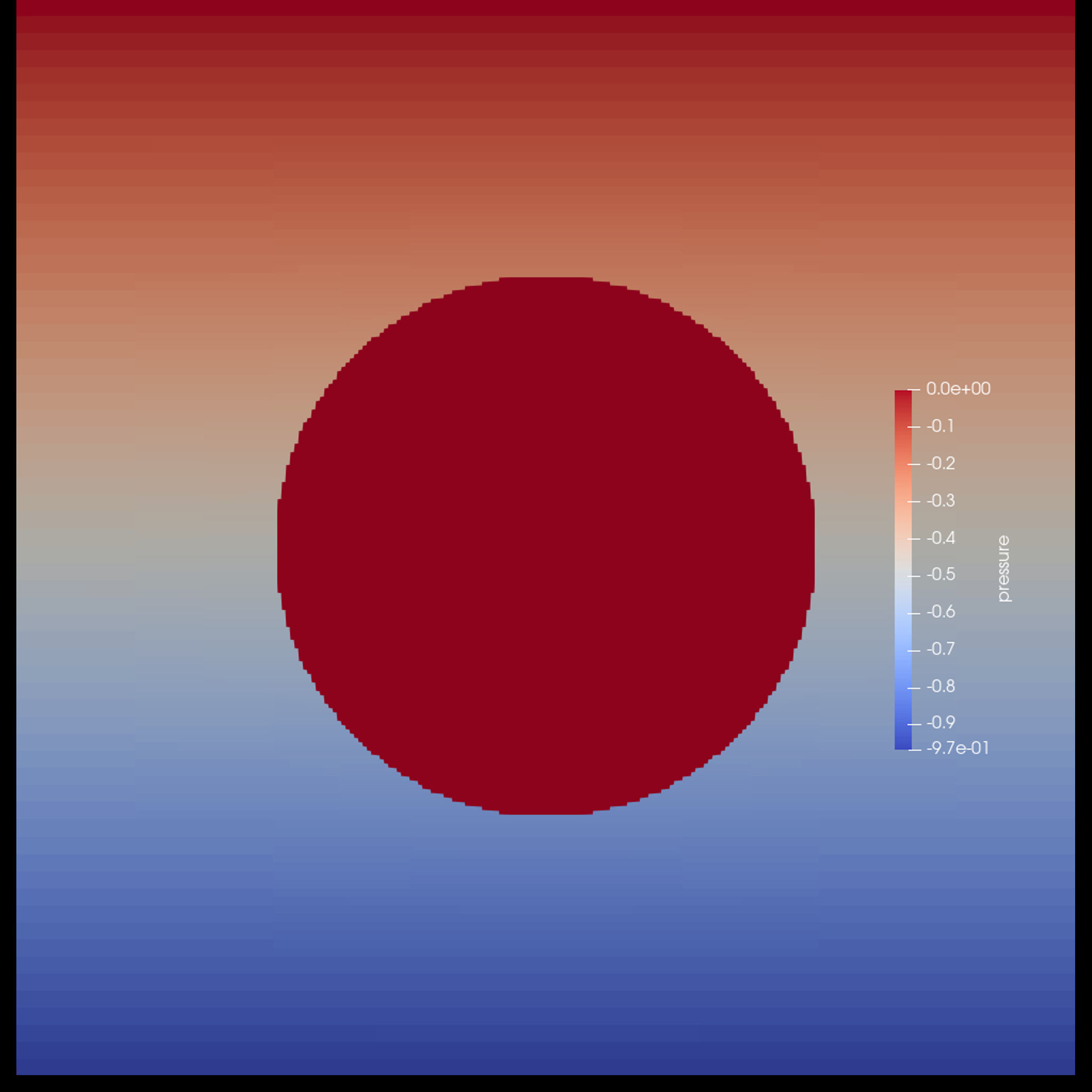}
	\end{minipage}%
	\hspace{0.02\textwidth}
	\begin{minipage}{0.38\textwidth}
		\centering
		\includegraphics[width=\linewidth]{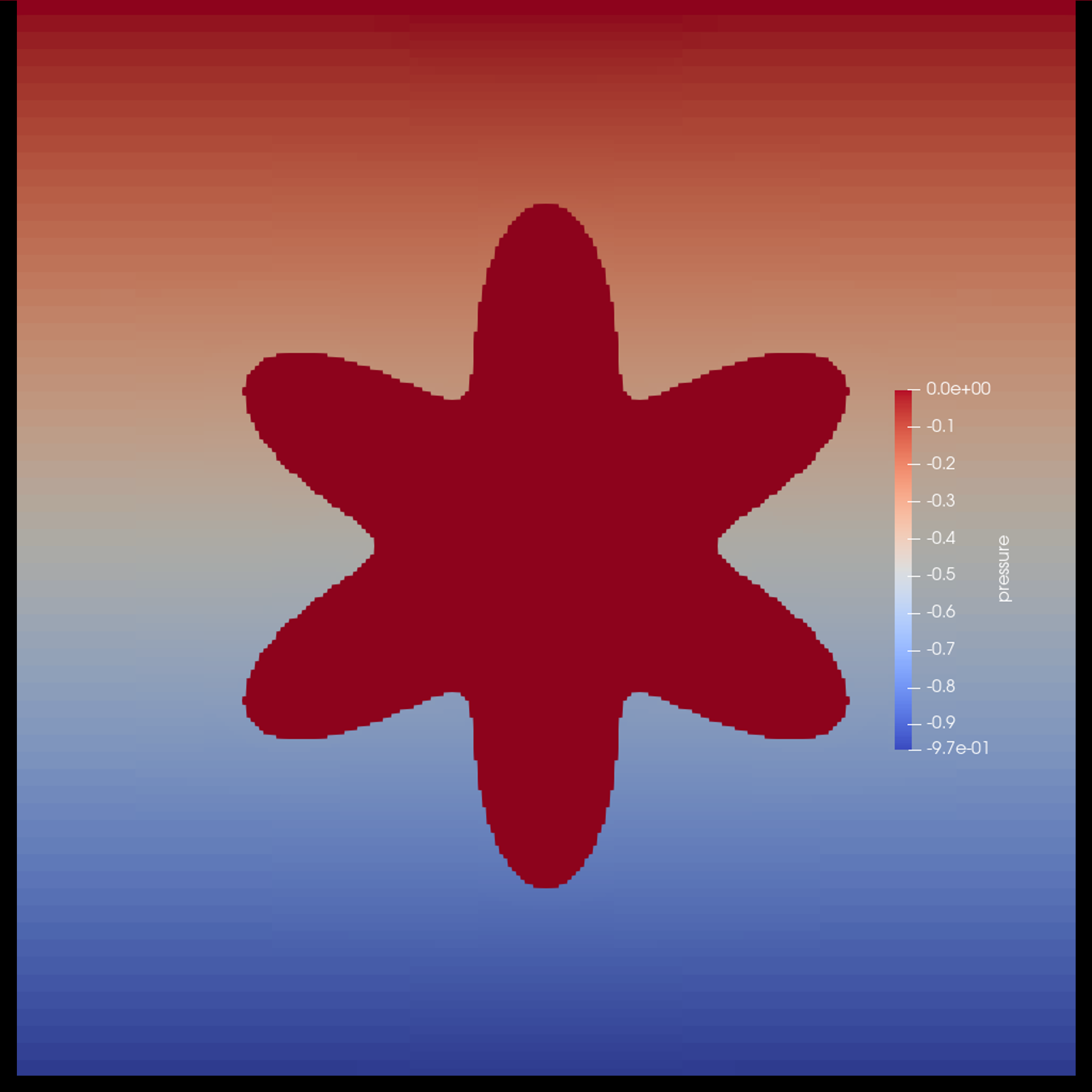}
	\end{minipage}%

	\caption{Pressure slice at $z=0.5$ for the static pressure projection test.
		The grids contain sphere (left) and star-shaped (right) obstacles with adaptive resolutions ranging from 64 to 256.
		Solid boundaries are shown in black.}
	\label{fig:cutcell_pressure}
\end{figure}

In this numerical test, we consider a scenario closer to practical fluid simulations.
A liquid tank occupying the computational domain $[0,1]^3$ is bounded by solid walls on the bottom and sides and open to air at the top.
A stationary solid obstacle of varying shape is placed at the center of the tank.
The purpose of this test is to evaluate the robustness of our pressure solver in the presence of complex solid boundaries.

The fluid is first assigned a uniform downward velocity $(0,-1,0)$.
We then solve the corresponding pressure Poisson system using multigrid-preconditioned PCG, where the multigrid hierarchy employs a $\mu$-cycle with $\mu=2$ to accelerate convergence.
The coefficients of the Poisson equation are constructed using the cut-cell method \cite{Ng2009efficient} to account for solid boundaries that may intersect grid cells arbitrarily.
After solving for pressure, we apply the pressure gradient to project the velocity field and measure the divergence of the projected velocity.

\autoref{fig:cutcell_pressure} shows slices of the pressure field computed by our solver for both the sphere and star-shaped obstacles.
The top boundary is open to air, and cells outside the liquid are treated as air by the solver, corresponding to a zero-pressure Dirichlet condition.
The resulting pressure field exhibits a nearly uniform gradient along the $y$-axis, which counteracts the imposed downward motion $(0,-1,0)$.

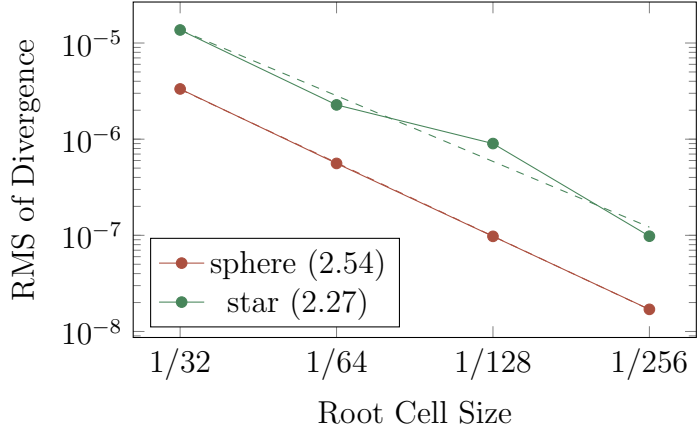
\begin{figure}[!htbp]
	\centering
	\begin{tikzpicture}
		\begin{semilogyaxis}[
				width=0.6\textwidth,
				height=0.4\textwidth,
				xlabel={Root Cell Size},
				ylabel=RMS of Divergence,
				xtick={1,2,3,4,5},
				xticklabels={
						1/16,
						1/32,
						1/64,
						1/128,
						1/256
					},
				cycle list={
						{color=myred, mark=*},
						{color=mygreen, mark=*},
						{color=myyellow, mark=*},
						{color=myblue, mark=*},
						{color=mypurple, mark=*},
						{color=mypink, mark=*}
					},
				yticklabel style={/pgf/number format/precision=0},
				legend pos=south west,
			]
			\addplot[color=myred, mark=*] coordinates {(2,3.3295531951539465e-06) (3,5.594628230774301e-07) (4,9.758516168553432e-08) (5,1.6984354445335158e-08)};
			\addlegendentry{sphere (2.54)};
			\addplot[domain=2:5, samples=100, dashed, thin, forget plot, color=myred] {2^(-13.138933616185556 + -2.536424375134655*x)};
			\addplot[color=mygreen, mark=*] coordinates {(2,1.3665405559265357e-05) (3,2.27258556020303e-06) (4,9.001816708929695e-07) (5,9.793827074868144e-08)};
			\addlegendentry{star (2.27)};
			\addplot[domain=2:5, samples=100, dashed, thin, forget plot, color=mygreen] {2^(-11.620016556006126 + -2.270936603038215*x)};
		\end{semilogyaxis}
	\end{tikzpicture}
	\caption{Volume-weighted RMS of divergence versus root cell size for the static pressure projection test.
		Dashed lines denote least-squares fitted convergence trends.
		The numbers in parentheses indicate the fitted grid convergence rates.}
	\label{fig:cutcell_resolution-wrms}
\end{figure}

\autoref{fig:cutcell_resolution-wrms} plots the volume-weighted RMS of divergence versus the root cell size.
Our solver achieves second-order convergence for all tested grids.

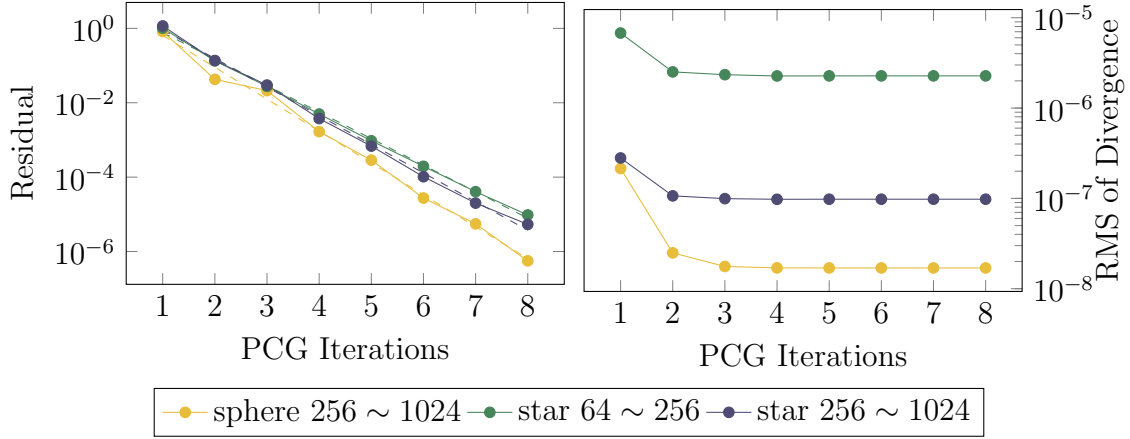
\begin{figure}[!htbp]
	\centering

	\begin{minipage}{0.49\linewidth}
		\centering
		\begin{tikzpicture}
			\begin{semilogyaxis}[
					width=\linewidth,
					height=0.72\linewidth,
					xlabel={PCG Iterations},
					ylabel={Residual},
					yticklabel style={/pgf/number format/precision=0},
					xtick={0,1,2,3,4,5,6,7,8},
					cycle list={
							{color=myred, mark=*},
							{color=mygreen, mark=*},
							{color=myyellow, mark=*},
							{color=myblue, mark=*},
							{color=mypurple, mark=*},
							{color=mypink, mark=*}
						},
					legend to name=staticpressurepcglegend,
					legend columns=4,
				]


				\addplot[color=myyellow, mark=*] coordinates {(1,8.236091e-01) (2,4.260168e-02) (3,2.116237e-02) (4,1.670045e-03) (5,2.859859e-04) (6,2.735510e-05) (7,5.543345e-06) (8,5.610053e-07)};
				\addlegendentry{sphere $256\sim 1024$};
				\addplot[domain=1:8, samples=100, dashed, thin, forget plot, color=myyellow] {2^(2.220915334899061 + -2.8484550610995245*x)};

				\addplot[color=mygreen, mark=*] coordinates {(1,1.014721e+00) (2,1.342742e-01) (3,2.751883e-02) (4,4.973482e-03) (5,9.503676e-04) (6,1.956842e-04) (7,4.062255e-05) (8,9.616197e-06)};
				\addlegendentry{star $64\sim 256$};
				\addplot[domain=1:8, samples=100, dashed, thin, forget plot, color=mygreen] {2^(1.9985252295569484 + -2.3697417611534926*x)};


				\addplot[color=mypurple, mark=*] coordinates {(1,1.162154e+00) (2,1.362108e-01) (3,2.989008e-02) (4,3.759314e-03) (5,6.830327e-04) (6,1.024869e-04) (7,2.000475e-05) (8,5.364198e-06)};
				\addlegendentry{star $256\sim 1024$};
				\addplot[domain=1:8, samples=100, dashed, thin, forget plot, color=mypurple] {2^(2.4222700828362402 + -2.5567357013001533*x)};

			\end{semilogyaxis}
		\end{tikzpicture}
	\end{minipage}
	\hfill
	\begin{minipage}{0.49\linewidth}
		\centering
		\begin{tikzpicture}
			\begin{semilogyaxis}[
					width=\linewidth,
					height=0.72\linewidth,
					xlabel={PCG Iterations},
					ylabel={RMS of Divergence},
					ytick pos=right,
					yticklabel pos=right,
					ylabel style={at={(axis description cs:1.20,0.5)},anchor=center},
					yticklabel style={/pgf/number format/precision=0},
					xtick={0,1,2,3,4,5,6,7,8},
					cycle list={
							{color=myred, mark=*},
							{color=mygreen, mark=*},
							{color=myyellow, mark=*},
							{color=myblue, mark=*},
							{color=mypurple, mark=*},
							{color=mypink, mark=*}
						},
				]


				\addplot[color=myyellow, mark=*] coordinates {(1,2.144079166401268e-07) (2,2.4949876059769416e-08) (3,1.7641471698672457e-08) (4,1.699967782716796e-08) (5,1.6982896324337648e-08) (6,1.6984183866478263e-08) (7,1.6984354445335158e-08) (8,1.6984383224803806e-08)};

				\addplot[color=mygreen, mark=*] coordinates {(1,6.769889505568633e-06) (2,2.515123635522988e-06) (3,2.34148548033982e-06) (4,2.2667996230787165e-06) (5,2.2709027384562705e-06) (6,2.2728227763803683e-06) (7,2.27265423072617e-06) (8,2.2725852809577508e-06)};


				\addplot[color=mypurple, mark=*] coordinates {(1,2.801713243887864e-07) (2,1.0665791837527716e-07) (3,9.937625733049109e-08) (4,9.775993134955878e-08) (5,9.791836450443217e-08) (6,9.794316271600427e-08) (7,9.79381686147054e-08) (8,9.793821992903118e-08)};

			\end{semilogyaxis}
		\end{tikzpicture}
	\end{minipage}

	\vspace{0.3em}

	\pgfplotslegendfromname{staticpressurepcglegend}

	\caption{PCG convergence for the static pressure projection test with cut cells on representative adaptive grids.
		Left: residual versus PCG iterations.
		Right: volume-weighted RMS of divergence versus PCG iterations.
		The divergence error reaches the discretization level after three iterations for all tested grids.
	}
	\label{fig:static-pressure-pcg-convergence}
\end{figure}

\autoref{fig:static-pressure-pcg-convergence} shows the convergence behavior of the multigrid-preconditioned PCG solver, terminated when the relative residual falls below $10^{-6}$.
Across all tested grids, the divergence error reaches the discretization level after approximately three PCG iterations.
Meanwhile, the residual decreases at a nearly grid-independent rate, with a per-iteration reduction factor larger than five.
The fitted reduction factors are $7.16$ for sphere $256\!\sim\!1024$, $5.17$ for star $64\!\sim\!256$, and $5.88$ for star $256\!\sim\!1024$.

\subsection{Performance and Solver Configuration}

\begin{table}[htbp]
	\caption{
		Performance of the solver on the sinusoidal Poisson test described in \autoref{subsec:sinusodial-poisson-test}.
		The numbers following each grid name indicate the range of refinement levels.
		All results are reported when the relative residual is reduced below $10^{-6}$.
		$T$ and $T_{\mathrm{iter}}$ denote the overall and per-iteration throughput, respectively.
		For small problem sizes, the GPU is not fully saturated and $T$ is correspondingly lower.
	}
	\label{tab:athena-sin-runtime-stats}
	\centering
	\begin{tabular}{lrrrr}
		\hline
		Grid          & Cells (M) & Time (ms) & $T$(Mcell/s) & $T_{\mathrm{iter}}$(Mcell/s) \\
		\hline
		sphere (3-5)  & 1.64453   & 6.795     & 242.03       & 1452.17                      \\
		star (3-5)    & 2.24609   & 8.109     & 276.98       & 1661.87                      \\
		uniform (4-4) & 2.00000   & 6.773     & 295.28       & 1771.69                      \\
		sphere (4-6)  & 7.55078   & 26.362    & 286.42       & 1718.53                      \\
		star (4-6)    & 10.96875  & 41.261    & 265.84       & 1595.04                      \\
		uniform (5-5) & 16.00000  & 69.566    & 230.00       & 1379.99                      \\
		sphere (5-7)  & 38.61328  & 179.436   & 215.19       & 1291.15                      \\
		star (5-7)    & 56.16797  & 266.064   & 211.11       & 1266.64                      \\
		\hline
	\end{tabular}
\end{table}

We evaluate the performance and configuration choices of our solver using two representative tests introduced earlier: the sinusoidal Poisson test and the static pressure projection test with cut cells.

First, we examine the full-solve throughput on the sinusoidal Poisson test.
\autoref{tab:athena-sin-runtime-stats} reports the total solve time and throughput $T$ for representative grid configurations, with the solver terminated when the relative residual falls below $10^{-6}$.
For large-scale problems, the full-solve throughput consistently exceeds $200$ million cells per second, demonstrating that the GPU implementation scales efficiently across different grid shapes and resolutions.

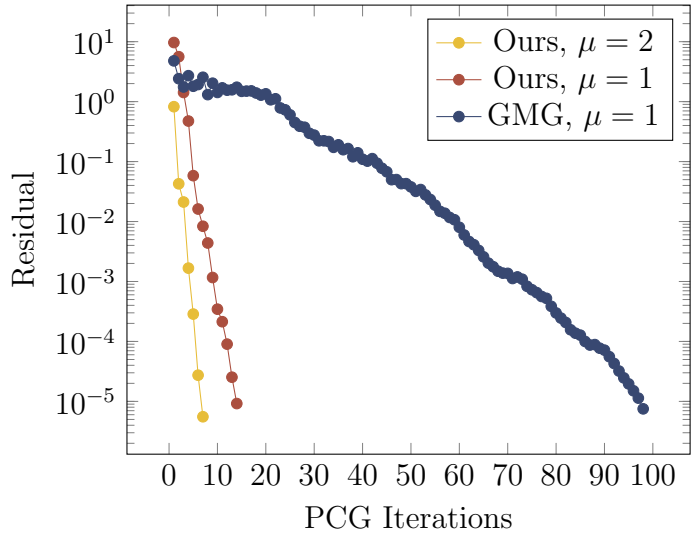
\begin{figure}[!htbp]
	\centering
	\begin{tikzpicture}
		\begin{semilogyaxis}[
				width=0.6\textwidth,
				height=0.5\textwidth,
				xlabel={PCG Iterations},
				ylabel={Residual},
				yticklabel style={/pgf/number format/precision=0},
				legend pos=north east,
				xtick={0,10, 20, 30, 40, 50, 60, 70,80,90,100},
				cycle list={
						{color=myred, mark=*},
						{color=mygreen, mark=*},
						{color=myyellow, mark=*},
						{color=myblue, mark=*},
						{color=mypurple, mark=*},
						{color=mypink, mark=*}
					},
			]

			\addplot[color=myyellow, mark=*] coordinates {(1,8.236091e-01) (2,4.260168e-02) (3,2.116237e-02) (4,1.670045e-03) (5,2.859859e-04) (6,2.735510e-05) (7,5.543345e-06)};
			\addlegendentry{Ours, $\mu=2$ };
			\addplot[color=myred, mark=*] coordinates {(1,9.680042e+00) (2,5.617069e+00) (3,1.421447e+00) (4, 4.743877e-01) (5, 5.825037e-02) (6, 1.616850e-02) (7, 8.340468e-03) (8, 4.385141e-03) (9, 1.165054e-03) (10, 3.454590e-04) (11, 2.133682e-04) (12, 9.018851e-05) (13, 2.532528e-05) (14, 9.201853e-06)};
			\addlegendentry{Ours, $\mu=1$ };
			\addplot[color=myblue, mark=*] coordinates {(1, 4.786161e+00) (2, 2.414754e+00) (3, 1.756567e+00)
					(4, 2.708138e+00) (5, 1.805904e+00) (6, 1.929030e+00) (7, 2.550695e+00)
					(8, 1.309044e+00) (9, 2.030225e+00) (10, 1.423748e+00) (11, 1.696602e+00)
					(12, 1.551985e+00) (13, 1.575682e+00) (14, 1.739187e+00) (15, 1.485600e+00)
					(16, 1.498159e+00) (17, 1.510439e+00) (18, 1.393031e+00) (19, 1.283093e+00)
					(20, 1.349449e+00) (21, 1.073733e+00) (22, 1.115700e+00) (23, 7.849178e-01)
					(24, 7.296770e-01) (25, 6.020103e-01) (26, 4.492221e-01) (27, 3.888050e-01)
					(28, 3.742207e-01) (29, 2.966285e-01) (30, 2.761435e-01) (31, 2.222453e-01)
					(32, 2.221427e-01) (33, 2.158929e-01) (34, 1.729110e-01) (35, 1.903938e-01)
					(36, 1.578002e-01) (37, 1.649304e-01) (38, 1.203662e-01) (39, 1.400573e-01)
					(40, 1.095507e-01) (41, 1.017718e-01) (42, 1.119862e-01) (43, 9.388231e-02)
					(44, 7.744215e-02) (45, 6.757937e-02) (46, 4.989358e-02) (47, 5.007952e-02)
					(48, 4.284567e-02) (49, 4.283199e-02) (50, 3.780477e-02) (51, 3.188070e-02)
					(52, 3.415344e-02) (53, 2.783854e-02) (54, 2.329673e-02) (55, 1.874312e-02)
					(56, 1.479449e-02) (57, 1.375334e-02) (58, 1.170515e-02) (59, 1.070386e-02)
					(60, 8.026253e-03) (61, 5.966288e-03) (62, 4.651166e-03) (63, 4.103838e-03)
					(64, 3.297497e-03) (65, 2.577297e-03) (66, 2.024583e-03) (67, 1.761910e-03)
					(68, 1.484790e-03) (69, 1.378306e-03) (70, 1.365714e-03) (71, 1.121785e-03)
					(72, 1.191695e-03) (73, 1.090893e-03) (74, 8.317118e-04) (75, 7.285895e-04)
					(76, 6.519382e-04) (77, 5.626399e-04) (78, 5.237872e-04) (79, 3.850726e-04)
					(80, 2.986614e-04) (81, 2.452488e-04) (82, 2.072272e-04) (83, 1.565239e-04)
					(84, 1.373260e-04) (85, 1.273719e-04) (86, 9.994500e-05) (87, 8.712287e-05)
					(88, 8.858709e-05) (89, 7.725472e-05) (90, 7.159241e-05) (91, 5.571158e-05)
					(92, 4.286628e-05) (93, 3.210324e-05) (94, 2.466906e-05) (95, 1.958793e-05)
					(96, 1.502244e-05) (97, 1.134205e-05) (98, 7.510436e-06)};
			\addlegendentry{GMG, $\mu=1$ };

		\end{semilogyaxis}
	\end{tikzpicture}
	\caption{Residual versus PCG iterations for the static pressure projection test with cut cells (star $256\sim1024$). The labels \emph{Ours, $\mu=1$} and \emph{Ours, $\mu=2$} correspond to our matrix-free multigrid preconditioner with algebraically consistent coarsening, using V-cycle and W-cycle respectively, while \emph{GMG} denotes geometric multigrid. Our method with the W-cycle ($\mu=2$) converges fastest, while GMG fails to converge effectively on this irregular grid.}
	\label{fig:static-pressure-pcg-residual-plot}
\end{figure}

Next, we study the effect of multigrid cycle configurations on the convergence behavior of the PCG solver using the static pressure projection test with cut cells.
\autoref{fig:static-pressure-pcg-residual-plot} compares several solver configurations.
Using our matrix-free multigrid preconditioner with a $\mu$-cycle ($\mu=2$) leads to the fastest convergence, reducing the residual by nearly two orders of magnitude within only a few PCG iterations.
Using a V-cycle ($\mu=1$) also converges, but requires more iterations.
In contrast, the geometric multigrid (GMG) preconditioner fails to converge effectively on this irregular cut-cell grid, leading to significantly slower residual reduction.








\subsection{Flow Around Kinematic Obstacles}
\label{sec:kinematic_flow}

To demonstrate the applicability of the proposed solver in practical fluid simulations with moving boundaries, we integrate the matrix-free multigrid-preconditioned PCG solver into the hybrid particle-grid flow map framework introduced in \citet{wang2025cirrus}. In this framework, a pressure Poisson equation must be solved at every timestep on adaptive octree grids.

The computational domain has size $1\times1\times2$ surrouded by Neumann boundary walls. A constant inflow velocity $\bm u=(0,0,1)$ is imposed at two walls $z=0$ and $z=2$, generating a vertical flow through the domain. The simulation runs for $2$ seconds at $200$ FPS, resulting in $400$ frames in total.

The adaptive grid hierarchy starts from a root grid of resolution $32\times32\times64$, and the finest level reaches an effective resolution of $512\times512\times1024$. The multigrid preconditioner uses a $\mu=2$ W-cycle. At each timestep, the Poisson equation defined by pressure projection is solved using the proposed matrix-free multigrid-preconditioned PCG solver.

We consider two kinematic obstacle configurations that follow a prescribed trajectory. The obstacle motion is defined as
\begin{equation}
	x = 0.5 + r\cos(2\pi s), \quad
	y = 0.5 + r\sin(2\pi s), \quad
	z = (1-s)z_0 + sz_1,
\end{equation}
where $s=t/2$ with $t\in[0,2]$, $r=0.2$, $z_0=1.7$, and $z_1=0.3$. This trajectory corresponds to a circular motion in the $x$--$y$ plane while the obstacle simultaneously descends along the $z$ direction. Over the two-second simulation interval, the obstacle completes one full revolution while moving downward through the domain.

In the first example, the obstacle is a sphere of radius $0.2$. The solid boundary is represented analytically using the signed distance function of the sphere. To estimate the fluid area fractions on cut-cell faces, the signed distance values are evaluated at cell corners and used to reconstruct the interface geometry.

In the second example, the analytic sphere is replaced by a triangle mesh representing a TIE fighter spacecraft model obtained from Sketchfab under CC BY license~\cite{sketchfab_tiefighter}. The mesh signed distance field is computed using the fast winding number method implemented in the libigl library. The obstacle follows the same translational trajectory as the sphere while simultaneously performing a rolling motion and maintaining a $10^\circ$ angle of attack. The complex geometry introduces significantly more irregular cut-cell configurations compared with the analytic sphere case.

\begin{figure}[htbp]
	\centering

	\begin{minipage}{0.38\textwidth}
		\centering
		\includegraphics[width=\linewidth]{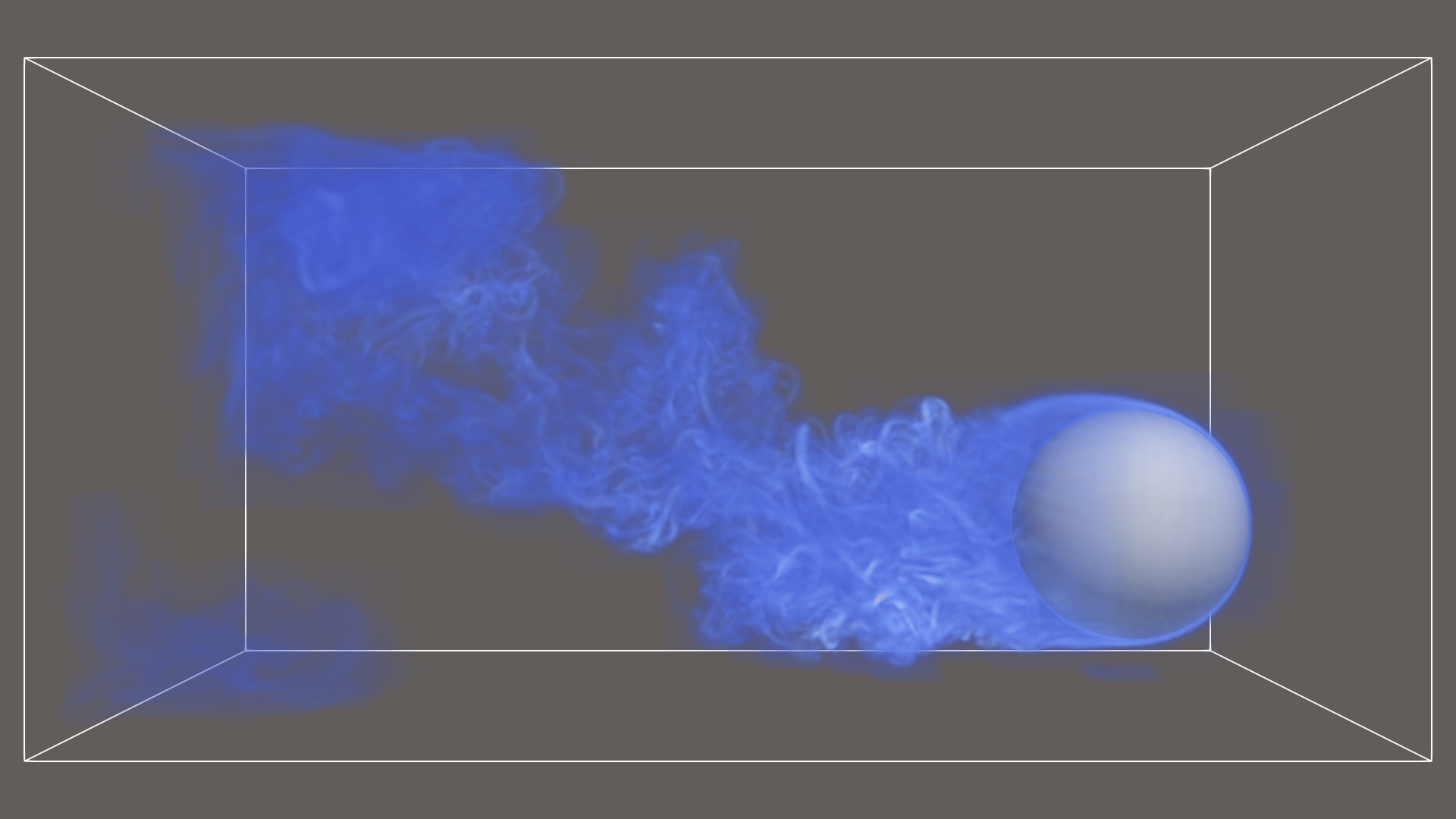}
	\end{minipage}%
	\hspace{0.02\textwidth}
	\begin{minipage}{0.38\textwidth}
		\centering
		\includegraphics[width=\linewidth]{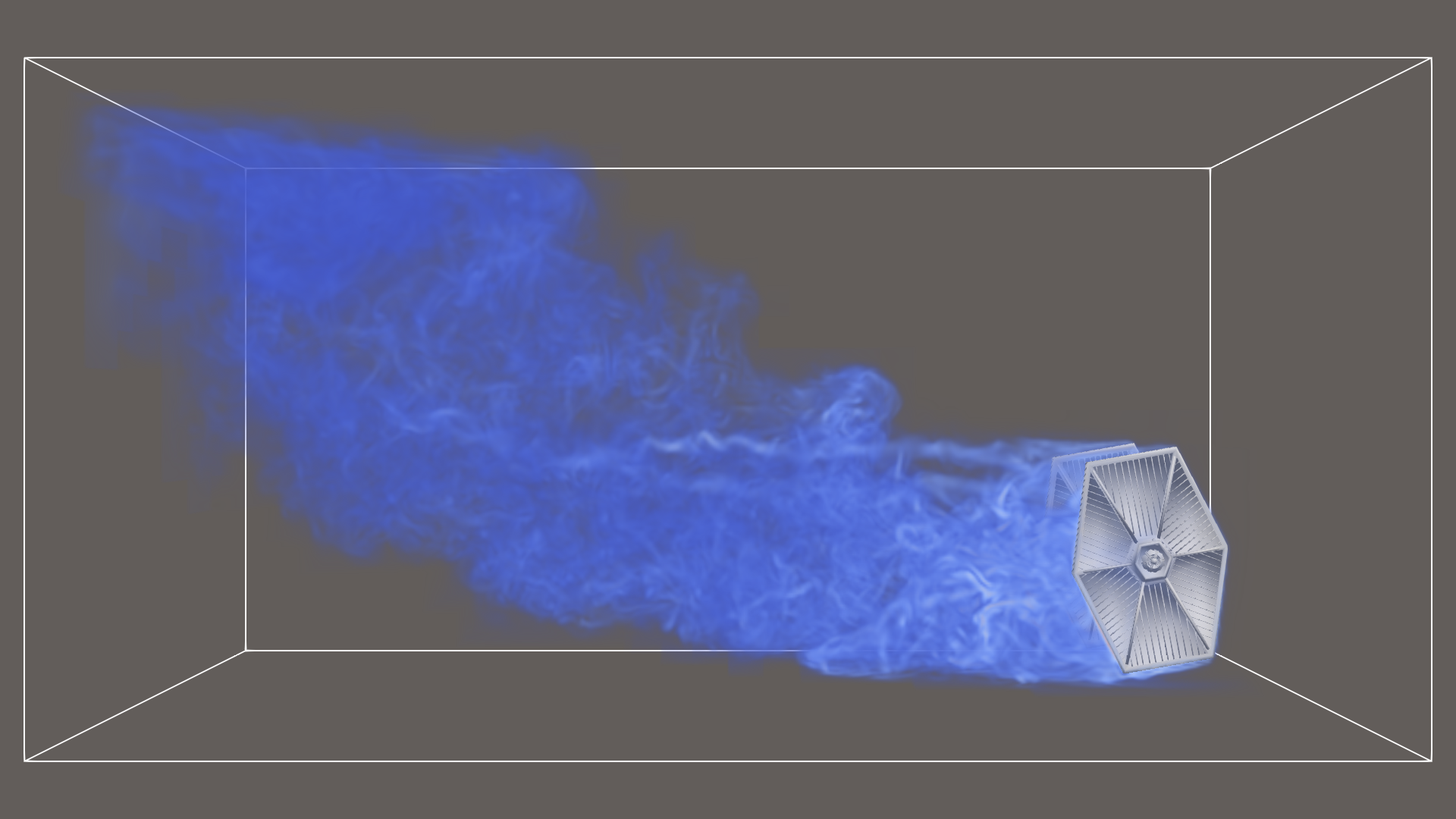}
	\end{minipage}%

	\caption{Volume rendering of the velocity magnitude at the final frame of the simulation.
		Left: flow around a kinematically moving sphere.
		Right: flow around a TIE fighter mesh following the same trajectory with additional rolling motion and a $10^\circ$ angle of attack.}

	\label{fig:kinematic_flow}
\end{figure}

\autoref{fig:kinematic_flow} shows volume renderings of the velocity magnitude at the final simulation frame for both examples. The results demonstrate that the solver remains stable and robust even in the presence of moving obstacles with complex geometries.

\autoref{tab:kinematic_flow_stats} reports the number of active leaf cells in the adaptive octree grid, the GPU wall time for the pressure projection step, the number of PCG iterations, and the full-solve throughput at the final timestep. In all cases, the pressure projection uses W-cycle ($\mu=2$) multigrid-preconditioned PCG, terminated when the relative residual falls below $10^{-6}$. The results indicate that the proposed solver maintains high computational efficiency despite the irregular domains introduced by moving obstacles.

For the relatively simple sphere case, the solver achieves a full-solve throughput of $\sim97$ M cell/s. The more geometrically complex TIE fighter case reaches $\sim72$ M cell/s, mainly due to the increased iteration count required by the more irregular cut-cell geometry. Nevertheless, as demonstrated in \autoref{subsec:sinusodial-poisson-test} and \autoref{subsec:static-projection-test}, the PCG solver typically reaches the discretization accuracy of the underlying numerical scheme within only about three iterations. In practical flow simulations, terminating the solver at this accuracy level could further reduce the runtime while maintaining the overall solution quality.


\begin{table}[htbp]
	\centering
	\caption{Performance statistics for the kinematic obstacle simulations. The reported time is the GPU wall time of the pressure projection step at the final timestep. $T$ denotes the overall solver throughput.}
	\label{tab:kinematic_flow_stats}
	\begin{tabular}{lcccc}
		\hline
		Case        & Leaf Cells & Time (ms) & PCG Iters & $T$(Mcell/s) \\
		\hline
		Sphere      & 11.61M     & 120       & 8         & 96.78        \\
		TIE fighter & 12.10M     & 169       & 11        & 71.62        \\
		\hline
	\end{tabular}
\end{table}

\section{Conclusion}
\label{sec:conclusion}

In this work, we presented a matrix-free multigrid preconditioner with
algebraically consistent coarsening for PCG applied to the Poisson equation on
adaptive octree grids with irregular domains. The method combines a compact
matrix-free discretization with a conservative cross-level formulation, and is
realized on the GPU through an efficient tile-based architecture.

The algebraically consistent coarsening plays a critical role in solving Poisson
equations on adaptive grids with cut cells. In contrast to purely geometric
multigrid methods, our approach enforces the Galerkin principle within
uniform-resolution patches, and introduces an FAS-style coarse-grid correction
at adaptive T-junctions to restore cross-level flux consistency where a strict
Galerkin construction is no longer available in the compact matrix-free form.
This strategy prevents the accumulation of additional residual errors across
refinement boundaries and coefficient discontinuities. As demonstrated in our
numerical experiments, the resulting PCG solver achieves grid-independent
convergence, typically reaching the discretization error level within only a few
iterations.

Furthermore, the matrix-free GPU implementation achieves high computational
efficiency. By storing only the diagonal coefficient and three face coefficients
per cell, the solver evaluates stencil operations on the fly and avoids the
storage and memory-access overhead of explicitly assembled sparse matrices. On a
single NVIDIA RTX 4090 GPU, the solver achieves full-solve throughputs above
200 million cells per second on analytical Poisson tests and above 70 million
cells per second on pressure projection problems in fluid simulation.

Despite its advantages, the proposed method has several limitations. First,
because the algebraically consistent coarsening is not available in the compact
matrix-free form at T-junctions, using a W-cycle ($\mu=2$) is often necessary to
achieve robust convergence in complex cut-cell configurations typical of fluid
simulations. Second, while the matrix-free coefficient storage is significantly
smaller than that required for an explicit sparse matrix, it still introduces
additional memory overhead compared with purely geometric multigrid formulations
where coefficients are recomputed on the fly.

Future work may explore more compact coefficient representations or alternative
transfer operators that better preserve the Galerkin structure across adaptive
interfaces while maintaining high GPU efficiency. Overall, the proposed approach
provides an efficient and practical Poisson solver for large-scale adaptive fluid
simulations and other applications involving elliptic equations on irregular
domains.



\bibliographystyle{ACM-Reference-Format}
\bibliography{cas-refs}





\end{document}